\documentclass[10pt,a4paper]{article}

\usepackage{amsmath,amsthm,amssymb}
\usepackage{calc}
\usepackage{bm} 

\newtheorem{thm}{\hspace*{20pt}Theorem}
\newtheorem{lem}{\hspace*{20pt}Lemma}[section]
\newtheorem{pro}{\hspace*{20pt}Proposition} 

\newcounter{constant}[section] 
\newcommand{\const}{\ifnum \theconstant >  0 \stepcounter{constant}\theconstant
                                             \else \setcounter{constant}{1}\theconstant \fi }
\newcounter{subconstant} 
\newcounter{submioneconstant} 
\newcounter{subconstanta}

\newcounter{submitwoconstant}

\newcommand{\suba}{\thesubconstanta} 
\newcommand{\insa}{\setcounter{subconstanta}{\value{constant}}}  
\newcounter{subconstantb} 
\newcommand{\subb}{\thesubconstantb} 
\newcommand{\insb}{\setcounter{subconstantb}{\value{constant}}}  
\newcounter{subconstantc} 
\newcommand{\subc}{\thesubconstantc} 
\newcommand{\insc}{\setcounter{subconstantc}{\value{constant}}}  
\newcounter{subconstantd} 
\newcommand{\subd}{\thesubconstantd} 
\newcommand{\insd}{\setcounter{subconstantd}{\value{constant}}}  
\newcounter{subconstante} 
\newcommand{\sube}{\thesubconstante} 
\newcommand{\inse}{\setcounter{subconstante}{\value{constant}}}  
\newcounter{subconstantf} 
\newcommand{\subf}{\thesubconstantf} 
\newcommand{\insf}{\setcounter{subconstantf}{\value{constant}}}  
\newcounter{subconstantg} 
\newcommand{\subg}{\thesubconstantg} 
\newcommand{\insg}{\setcounter{subconstantg}{\value{constant}}}  
\newcounter{subconstanth} 
\newcommand{\subh}{\thesubconstanth} 
\newcommand{\insh}{\setcounter{subconstanth}{\value{constant}}}  
\newcounter{subconstanti} 
\newcommand{\subi}{\thesubconstanti} 
\newcommand{\insi}{\setcounter{subconstanti}{\value{constant}}}  
\newcounter{subconstantj} 
\newcommand{\subj}{\thesubconstantj} 
\newcommand{\insj}{\setcounter{subconstantj}{\value{constant}}}  
\newcounter{subconstantk} 
\newcommand{\subk}{\thesubconstantk} 
\newcommand{\insk}{\setcounter{subconstantk}{\value{constant}}}  
\newcounter{subconstantl}

\newcounter{subconstantm}

\newcounter{subconstantn}

\newcounter{subconstanto}

\newcounter{subconstantp}

\newcounter{subconstantq}

\newcounter{subconstantr}

\newcounter{subconstants}

\newcounter{subconstantt}

\newcounter{subconstantu}

\newcounter{subconstantv}

\newcounter{subconstantw}

\newcounter{subconstantx}

\newcounter{subconstanty} 
\newcommand{\suby}{\thesubconstanty} 
\newcommand{\insy}{\setcounter{subconstanty}{\value{constant}}}   
\newcounter{subconstantz} 
\newcommand{\subz}{\thesubconstantz} 
\newcommand{\insz}{\setcounter{subconstantz}{\value{constant}}}

\makeatletter
\@addtoreset{equation}{section}

\makeatother

\begin{document}

\title{Critical blowup exponent to a class of semilinear elliptic equations with constraints in higher dimension - local properties}




\author{Takashi Suzuki\footnotemark[1] and Ryo Takahashi \footnotemark[2]}

\footnotetext[1]{
Division of Mathematical Science, Department of Systems Innovation, 
Graduate School of Engineering Science, Osaka University, 
Machikaneyamacho 1-3, Toyonakashi, 560-8531, Japan. 
(E-mail: {\it suzuki@sigmath.es.osaka-u.ac.jp})
}
\footnotetext[2]{
Division of Mathematical Science, Department of Systems Innovation, 
Graduate School of Engineering Science, Osaka University, 
Machikaneyamacho 1-3, Toyonakashi, 560-8531, Japan. 
(E-mail: {\it r-takaha@sigmath.es.osaka-u.ac.jp})
}

\date{\today}

\maketitle

\begin{abstract}
We study a class of semilinear elliptic equations with constraints in higher dimension. 
It is known that several mathematical structures of the problem are closed to those of the Liouville equation in dimension two. 
In this paper, we establish a classification of entire solutions, the $\sup + \inf$ type inequality and the quantized blowup mechanism. 
\end{abstract}


\section{Introduction}
In this paper, we are concerned with the semilinear elliptic equation 
\begin{equation}
\begin{cases}
& -\Delta v = A(x)v_+^\gamma \quad \mbox{in $\Omega$} \\ 
& \int_{\Omega} v_+^{\frac{n(\gamma -1)}{2}} dx \leq T < +\infty, 
\end{cases}
 \label{eqn:higher-d}
\end{equation}
where $v=v(x)$ is an unknown function, $\gamma \in \left(1, \frac{n+2}{n-2} \right)$, $n \geq 3$, $v_+=\max\{v,0\}$, $\Omega \subset {\bf R}^n$ is an open set, 
$T\geq 0$ is a given constant, and $A=A(x)$ is a function defined on $\overline{\Omega}$. 
We are also interested in the problem
\begin{equation}
\begin{cases}
& -\Delta v = v_+^\gamma \quad \mbox{in ${\bf R}^n$} \\ 
& \int_{\Omega} v_+^{\frac{n(\gamma -1)}{2}} dx \leq T < +\infty. 
\end{cases}
 \label{eqn:higher-d'}
\end{equation}
Equation (\ref{eqn:higher-d}) arises in various situations. 
For example, (\ref{eqn:higher-d}) appears as a free boundary value problem when plasma confinement is considered, see \cite{te75} and \cite{te77}. 
Also, in astrophysics, the degenerate parabolic equation is derived from the kinetic theory, see \cite{chava} and \cite{cs04}, 
and there (\ref{eqn:higher-d}) is derived from the total mass conservation and the decrease of the free energy, see \cite{suzuki07} and \cite{st2}. 

As pointed out in \cite{wy03}, equation (\ref{eqn:higher-d}) has several properties similar to the Liouville equation in dimension two, which is described by 
\begin{equation}
\begin{cases}
& -\Delta u = V(x)e^u \quad \mbox{in $\Omega\subset{\bf R}^2$} \\ 
& \int_{\Omega} e^u dx < +\infty, 
\end{cases}
 \label{eqn:liouville}
\end{equation}
where $u=u(x)$ is an unknown function, $\Omega\subset {\bf R}^2$ is an open set, and $V=V(x)$ is a function defined on $\overline{\Omega}$. 

At first, we notice that both of (\ref{eqn:higher-d}) and (\ref{eqn:liouville}) are invariant under the scalings
\[
\mu^q v(\mu x) \quad \mbox{and} \quad u(\mu x)+2\log \mu 
\]
for $\mu>0$, respectively, where $q=\frac{2}{\gamma-1}$. 
By virtue of these scaling invariance, we can develop the blowup analysis for both of (\ref{eqn:higher-d}) and (\ref{eqn:liouville}). 
It is known that the blowup analysis works well for (\ref{eqn:liouville}). 
We find later that this situation is also applicable to (\ref{eqn:higher-d}). 

To develop the blowup analysis, it is essential to classify the entire solutions. 
As to (\ref{eqn:liouville}) for $\Omega={\bf R}^2$ and $V\equiv 1$, Chen and Li (\cite{cl91}) showed that any nontrivial solution is explicitly given by
\[
v(x) = \log\left\{ \frac{8\mu^2}{(1+\mu^2 |x-x_0|^2)^2} \right\}
\]
for some $x_0\in{\bf R}^2$ and $\mu>0$. 
Also, Wang and Ye (\cite{wy03}) classified nontrivial solutions of (\ref{eqn:higher-d'}) for $\gamma=\frac{n}{n-2}$. 
The first aim of this paper is to extend this result of \cite{wy03} to the case that $\gamma \in \left(1, \frac{n+2}{n-2} \right)$. 
To state the first result, we introduce several notations. 
Let $\phi=\phi(r)$ be the unique classical solution to the problem 
\begin{equation}
\begin{cases}
& \phi''+\frac{n-1}{r}\phi'+\phi_+^\gamma=0,\ r>0, \\ 
& \phi(0)=1,\quad \phi'(0)=0, 
\end{cases}
 \label{eqn:phi-ode}
\end{equation}
and let $r_\gamma^\ast>0$ be the first zero point of $\phi=\phi(r)$. 
Then, there exists $\alpha_\gamma^\ast>0$ such that $\phi'(r_\gamma^\ast)=-\alpha_\gamma^\ast<0$. 
We put 
\begin{equation}
\lambda_\gamma^\ast=\omega_{n-1}\int_0^{r_\gamma^\ast}\phi^\gamma r^{n-1}dr=\omega_{n-1}\alpha_\gamma^\ast(r_\gamma^\ast)^{n-1}, 
 \label{eqn:def-M}
\end{equation}
where $\omega_{n-1}$ stands for the area of the boundary of the unit ball in ${\bf R}^n$. 
Under these preparations, the first result is stated as follows. 

\begin{thm}\label{thm:class}
Assume that $\gamma \in \left(1, \frac{n+2}{n-2} \right)$ and $n \geq 3$. 
Then, any nontrivial classical solution of (\ref{eqn:higher-d'}) is radially symmetric about some point, and satisfies 
\[
\int_{{\bf R}^n} v_+^{\frac{n(\gamma-1)}{2}} dx=\lambda_\gamma^\ast. 
\]
More precisely, $v=v(x)$ is represented by 
\begin{equation}
v(x)=\begin{cases}
\mu^q\phi(\mu|x-x_0|) & (|x-x_0|\leq r_\gamma^\ast/\mu) \\
\frac{\mu^{q-(n-2)}\lambda_\gamma^\ast}{\omega_{n-1}(n-2)}\left(\frac{1}{|x-x_0|^{n-2}}-\frac{1}{(r_\gamma^\ast/\mu)^{n-2}}\right) & (|x-x_0|\geq r_\gamma^\ast/\mu)
\end{cases}
 \label{eqn:thm:class-1}
\end{equation}
for some $x_0\in{\bf R}^n$ and $\mu>0$, where $q=\frac{2}{\gamma-1}$. 
\end{thm}

Let $\{u_k\}$ be a solution sequence of (\ref{eqn:liouville}) for $V=V_k\in C(\overline{\Omega})$, 
and assume that there exists $C>0$ and $V\in C(\overline{\Omega})$ such that 
\[
\int_{\Omega} e^{u_k}dx\leq C,\quad V_k\geq 0\ \mbox{in $\Omega$},\quad V_k\rightarrow V \ \mbox{in $C(\overline{\Omega})$}. 
\]
Then, passing to a subsequence, we have the following alternatives: 

{\bf (i)} $\{u_k\}$ is locally uniformly bounded in $\Omega$. 

{\bf (ii)} $u_k \rightarrow -\infty$ locally uniformly in $\Omega$. 

{\bf (iii)} There exist $l$-points $\{x_i\}_{i=1}^l$ such that 
$v_k \rightarrow -\infty$ locally uniformly in $\Omega\setminus\{x_1,\cdots,x_l\}$ and  
\[
V_k(x)e^{u_k}dx \overset{\ast}{\rightharpoonup} \sum_{i=1}^{l} \alpha(x_i) \delta_{x_i}(dx) \quad \mbox{in ${\cal M}(\Omega)$}
\]
with $\alpha(x_i)\in 8\pi{\bf N}$ for $i=1,\cdots,l$, 
where $\delta_x$ and ${\cal M}(\Omega)$ denote the Dirac measure centered at $x$ and the space of measure identified with the dual space of $C(\Omega)$, respectively. 

This property, called {\it quantized blowup mechanism} in this paper, is proven in \cite{bm91}, except for the sharp result $\alpha(x_i)\in 8\pi{\bf N}$ which is proven in \cite{ls94}. 
An obstacle to prove the sharp result is to show the {\it residual vanishing}. 
The key to show it is the $\sup+\inf$ type inequality shown in \cite{shaf92} and \cite{bls93}. 
Inequality of this type to (\ref{eqn:higher-d}) is shown for the case that $\gamma=\frac{n}{n-2}$ and $A\equiv 1$ in \cite{wy03}. 
We here extend it to the case that $\gamma \in \left(1, \frac{n+2}{n-2} \right)$ with a perturbation $A=A(x)$. 

\begin{thm}\label{thm:supinf}
Assume that $\gamma \in \left(1, \frac{n+2}{n-2} \right)$, $n \geq 3$, $\Omega$ is an open set, $A\in C(\overline{\Omega})$, and $0\leq A\leq C_{\const\insc}$ in $\Omega$ for some $C_{\subc}>0$. 
Then, for any compact set $K \subset \Omega$ and any number $T>0$, 
there exist $C_{\const\insa}=C_{\suba}(n,\gamma,A)>0$ and $C_{\const\insb}=C_{\subb}(n,\gamma,K,T,A)>0$ 
such that 
\[
\sup_K v+C_{\suba}\inf_\Omega v\leq C_{\subb}
\]
for any solution $v=v(x)$ of (\ref{eqn:higher-d}). 
\end{thm}

Quantized blowup mechanism to (\ref{eqn:higher-d}) is shown for the case that $\gamma=\frac{n}{n-2}$ and $A\equiv 1$ in \cite{wy03}. 
From the blowup analysis based on Theorem \ref{thm:class}, 
we obtain the quantized blowup mechanism even for the case that $\gamma \in \left(1, \frac{n+2}{n-2} \right)$ with perturbation $A=A(x)$. 

\begin{thm}\label{thm:quantization} 
Assume that $\gamma \in \left(1, \frac{n+2}{n-2} \right)$, $n \geq 3$, $\Omega$ is an open set, 
$0\leq A_k\in C(\overline{\Omega})$, 
and there exists $0\leq A\in C(\overline{\Omega})$ such that $A_k\rightarrow A$ in $C(\overline{\Omega})$. 
Then, given $T>0$ and a solution sequence $v_k=v_k(x)$ of 
\begin{equation}
\begin{cases}
-\Delta v_k = A_k(x)(v_k)_+^\gamma \quad \mbox{in $\Omega$} \\ 
\int_{\Omega} (v_k)_+^{\frac{n(\gamma -1)}{2}} dx \leq T,  
\end{cases}
 \label{eqn:bm-ls-type}
\end{equation}
passing to a subsequence, we have the following alternatives: 

{\bf (i)} $\{ v_k \}$ is locally uniformly bounded in $\Omega$. 

{\bf (ii)} $v_k \rightarrow -\infty$ locally uniformly in $\Omega$. 

{\bf (iii)} There exists a finite set ${\cal S} = \{x_i\}_{i=1}^l$ such that 
$v_k \rightarrow -\infty$ locally uniformly in $\Omega\setminus{\cal S}$, and that 
\[
(v_k)_+^{\frac{n(\gamma -1)}{2}}dx \overset{\ast}{\rightharpoonup} \sum_{i=1}^{l} m(x_i) \delta_{x_i}(dx) \quad \mbox{in ${\cal M}(\Omega)$}
\]
with $m(x_i)\in A(x_i)^{-n/2}\lambda_\gamma^\ast{\bf N}$ for $i=1,\cdots,l$, 
where $\lambda_\gamma^\ast$ is as in Theorem \ref{thm:class}. 
\end{thm}

From parallel properties between (\ref{eqn:higher-d}) and (\ref{eqn:liouville}), 
it is expected that some advanced results following the theorems stated above hold, see \cite{suzuki07} and the references therein for advanced results corresponding to (\ref{eqn:liouville}). 
Actually, we study the asymptotic profile of solution sequences and the location of the blowup points to equation (\ref{eqn:higher-d}) with free boundary value in the forthcoming paper, 
which corresponds to the results of \cite{ns90} and \cite{st} for equation (\ref{eqn:liouville}). \\

This paper is composed of five sections. 
Some preparatory lemmas are provided and shown in Section \ref{sec:pre}, 
and then Theorem \ref{thm:class} is proven in Section \ref{sec:class}. 
We prove Theorems \ref{thm:supinf} and \ref{thm:quantization} in Sections \ref{sec:supinf} and \ref{sec:bm-ls-type}, respectively. 

Henceforth, $C_i$ ($i=1,2,\cdots$) denote positive constants whose subscripts are renewed in each section, 
and we shall not distinguish any sequences with their subsequences for the sake of shorthand. 
\section{Preliminaries}\label{sec:pre}
In this section, we shall provide some preparatory lemmas to prove the theorems stated in the previous section. 

We shall use the following lemma to show the radial symmetry of Theorem \ref{thm:class}. 
Although the lemma is a part of the result of \cite{li-n93}, we here give the proof for completeness. 
\begin{lem}\label{lem:pre-1}
For every $\gamma>1$ and $\sigma>0$, any classical solution of 
\[
\begin{cases}
-\Delta w=(w-\sigma)_+^\gamma,\ w>0 & \mbox{in ${\bf R}^n$} \\
w(x)\rightarrow 0 & \mbox{as $|x|\rightarrow+\infty$}. 
\end{cases}
\]
is radially symmetric about some $x_0\in{\bf R}^n$ and satisfies $\partial u/\partial r<0$ for $r=|x-x_0|>0$. 
\end{lem}

{\it Proof.}\ 
Given $\gamma>1$ and $\sigma>0$, let $w=w(x)$ be a classical solution of 
\[
\begin{cases}
-\Delta w=(w-\sigma)_+^\gamma,\ w>0 & \mbox{in ${\bf R}^n$} \\
w(x)\rightarrow 0 & \mbox{as $|x|\rightarrow+\infty$}. 
\end{cases}
\]
For $x=(x_1,x_2,\cdots,x_n)\in{\bf R}^n$ and $\lambda\in{\bf R}$, we define
\[
T_\lambda=\{x\in{\bf R}^n \ | \ x_1=\lambda\},\quad 
\Sigma_\lambda=\{x\in{\bf R}^n \ | \ x_1<\lambda\},\quad 
x^\lambda=(2\lambda-x_1,x_2,\cdots,x_n). 
\]
We also introduce 
\[
\Lambda=\{\lambda\in{\bf R} \ | \ w(x)>w(x^\lambda) \ \mbox{for $x\in\Sigma_\lambda$ and } \partial w/\partial x_1<0 \ \mbox{on $T_\lambda$}\}. 
\]
Since $w>0$ in ${\bf R}^n$ and $w(x)\rightarrow 0$ as $|x|\rightarrow+\infty$, 
there exists $0<R_0<R_1$ such that 
\begin{align}
& \max_{{\bf R}^n\setminus B_{R_0}} w\leq \frac{\sigma}{2}
 \label{eqn:pre-1-1}\\
& 0<\max_{{\bf R}^n\setminus B_{R_1}} w\leq \frac{1}{2}\min_{\overline{B_{R_0}}} w. 
 \label{eqn:pre-1-2}
\end{align}
We now show the following properties. \\

(i) $[R_1,+\infty)\subset\Lambda$. 

(ii) For any $\lambda_0\in\Lambda\cap(0,+\infty)$, there exists $\varepsilon_0>0$ such that 
$(\lambda_0-\varepsilon_0,\lambda_0+\varepsilon_0)\subset\Lambda\cap(0,+\infty)$. 

(iii) We have either $w(x)\equiv w(x^{\lambda_1})$ in $\Sigma_{\lambda_1}$ for some $\lambda_1\geq 0$, or
\[
w(x_1,x')>w(-x_1,x')\ \mbox{for $x_1<0$},\quad 
\frac{\partial w}{\partial x_1}(x_1,x')<0 \ \mbox{for $x_1>0$}, 
\]
where $x'=(x_2,\cdots,x_n)$. \\

[{\it Proof of (i)}] Fix $\lambda\geq R_1$ and put $z(x)=w(x)-w(x^\lambda)$. 
Then we have
\begin{equation}
z>0 \quad \mbox{in $\overline{B_{R_0}}$}
 \label{eqn:pre-1-3}
\end{equation}
since it holds by (\ref{eqn:pre-1-2}) that 
\[
z(x)=w(x)-w(x^\lambda)\geq \min_{\overline{B_{R_0}}}w-\max_{{\bf R}^n\setminus B_{R_1}}w\geq \frac{1}{2}\min_{\overline{B_{R_0}}}w>0
\]
for $x\in\overline{B_{R_0}}$. 
Moreover, it follows from (\ref{eqn:pre-1-1}) and (\ref{eqn:pre-1-3}) that 
\[
\begin{cases}
-\Delta z=0 \ \mbox{in $\Sigma_\lambda\setminus\overline{B_{R_0}}$} \quad 
z\geq 0 \ \mbox{on $\partial(\Sigma_\lambda\setminus\overline{B_{R_0}})$} \\
z(x)\rightarrow 0 \ \mbox{as $|x|\rightarrow+\infty$, $x\in\Sigma_\lambda\setminus\overline{B_{R_0}}$}. 
\end{cases}
\]
Since $z\not\equiv 0$ in $\Sigma_\lambda\setminus\overline{B_{R_0}}$, the maximum principle and the Hopf lemma assure 
\begin{equation}
z>0 \ \mbox{in $\Sigma_\lambda\setminus\overline{B_{R_0}}$}, \quad 
\frac{\partial z}{\partial x_1}=2\frac{\partial w}{\partial x_1}<0 \ \mbox{on $T_\lambda$}. 
 \label{eqn:pre-1-4}
\end{equation}
Thus (\ref{eqn:pre-1-3}) and (\ref{eqn:pre-1-4}) yield property (i). \\

[{\it Proof of (ii)}] We may assume $\lambda_0\leq R_1$ by (i). 
Since
\[
w(x)-w(x^{\lambda_0})>0 \ \mbox{for $x\in \Sigma_{\lambda_0}$}, \quad 
\frac{\partial w}{\partial x_1}<0 \ \mbox{on $T_{\lambda_0}$}, 
\]
we see from the continuity of $\partial w/\partial x_1$ that there exists $0<\varepsilon_1\ll 1$ such that 
\[
\frac{\partial w}{\partial x_1}<0 \ \mbox{in}\ \{x=(x_1,x')\in{\bf R}^n \ | \ x_1\in[\lambda_0-4\varepsilon_1,\lambda_0+4\varepsilon_1],\ |x'|\leq R_1+1\}, 
\]
which implies 
\begin{equation}
\begin{cases}
w(x)-w(x^\lambda)>0 & \mbox{in } \{x\in\overline{B_{R_1+1}} \ | \ \lambda_0-2\varepsilon_1\leq x_1<\lambda\} \\
\frac{\partial w}{\partial x_1}<0 & \mbox{on } \overline{B_{R_1+1}}\cap T_\lambda
\end{cases}
 \label{eqn:pre-1-5}
\end{equation}
for any $\lambda\in(\lambda_0-\varepsilon_1,\lambda_0+\varepsilon_1)$. 
In addition, for 
\begin{align*}
&M=2\max\{|\partial w/\partial x_1| \ | \ |x_1|\leq 2(R_1+1),\ |x'|\leq R_1+1\}>0 \\
&\delta=\min\{w(x)-w(x^{\lambda_0}) \ | \ |x'|\leq R_1+1,\ x_1\in[-(R_1+1),\lambda_0-2\varepsilon_1]\}>0, 
\end{align*}
we find
\begin{equation}
w(x)-w(x^\lambda)\geq \frac{\delta}{2}>0 \quad \mbox{in } \{x\in\overline{B_{R_1+1}} \ | \ x_1\in[-(R_1+1),\lambda_0-2\varepsilon_1]\}
 \label{eqn:pre-1-6}
\end{equation}
for any $\lambda\in(\lambda_0-\varepsilon_0,\lambda_0+\varepsilon_0)$, where 
\[
\varepsilon_0=\min\left\{\varepsilon_1,\frac{\delta}{2M},\lambda_0\right\}. 
\] 
Combining (\ref{eqn:pre-1-5}) and (\ref{eqn:pre-1-6}) shows
\begin{equation}
w(x)-w(x^\lambda)>0 \ \mbox{in $\overline{B_{R_1+1}}\cap\Sigma_\lambda$},\quad 
\frac{\partial w}{\partial x_1}<0 \ \mbox{on $\overline{B_{R_1+1}}\cap T_\lambda$}
 \label{eqn:pre-1-7}
\end{equation}
for any $\lambda\in(\lambda_0-\varepsilon_0,\lambda_0+\varepsilon_0)$. 
On the other hand, for every $\lambda\in(\lambda_0-\varepsilon_0,\lambda_0+\varepsilon_0)$, $z(x)=w(x)-w(x^\lambda)$ satisfies 
\[
\begin{cases}
-\Delta z=0,\ z\neq 0\ \mbox{in $\Sigma_\lambda\setminus\overline{B_{R_1+1}}$},\quad 
z\geq 0\ \mbox{on $\partial(\Sigma_\lambda\setminus\overline{B_{R_1+1}})$} \\
z(x)\rightarrow 0 \ \mbox{as $|x|\rightarrow+\infty$, $x\in\Sigma_\lambda\setminus\overline{B_{R_1+1}}$}
\end{cases}
\]
by (\ref{eqn:pre-1-1}) and (\ref{eqn:pre-1-7}), which implies 
\begin{equation}
w(x)-w(x^\lambda)>0 \ \mbox{in $\Sigma_\lambda\setminus\overline{B_{R_1+1}}$},\quad 
\frac{\partial w}{\partial x_1}<0 \ \mbox{on $T_\lambda\setminus\overline{B_{R_1+1}}$}
 \label{eqn:pre-1-8}
\end{equation}
for any $\lambda\in(\lambda_0-\varepsilon_0,\lambda_0+\varepsilon_0)$. 
Claim (ii) follows from (\ref{eqn:pre-1-7}) and (\ref{eqn:pre-1-8}). \\

[{\it Proof of (iii)}] We put 
\[
\lambda_\ast=\inf\{\lambda\geq 0 \ | \ (\lambda,+\infty)\subset\Lambda \}. 
\]
By the continuity of $w$ and the definition of $\Lambda$, we have 
\[
z(x)=w(x)-w(x^{\lambda_\ast})\geq 0 \quad \mbox{in $\Sigma_{\lambda_\ast}$}. 
\]
It follows from the mean value theorem that there exists $c=c(x)\in C(\overline{\Sigma_{\lambda_\ast}})$ such that 
\[
\begin{cases}
-\Delta z=cz,\ z\geq 0\ \mbox{in $\Sigma_{\lambda_\ast}$},\quad 
z=0\ \mbox{on $T_{\lambda_\ast}$} \\
z(x)\rightarrow 0 \ \mbox{as $|x|\rightarrow+\infty$, $x\in\Sigma_{\lambda_\ast}$}. 
\end{cases}
\]
Hence it holds either 
\[
w(x)=w(x^{\lambda_\ast}) \quad \mbox{in $\Sigma_{\lambda_\ast}$}
\]
or 
\[
w(x)-w(x^{\lambda_\ast})>0\ \mbox{in $\Sigma_{\lambda_\ast}$},\quad 
\frac{\partial w}{\partial x_1}<0 \ \mbox{on $T_{\lambda_\ast}$}
\]
by the maximum principle and the Hopf lemma. 
If the latter and $\lambda_\ast>0$ simultaneously occur, 
then $(\lambda_\ast-\varepsilon,+\infty)\subset\Lambda$ for some $0<\varepsilon<\lambda_\ast$ by (ii). 
But this contradicts the definition of $\lambda_\ast$, and therefore $\lambda_\ast=0$ when the latter holds. \\ 

We are now in a position to prove the lemma. 
If $w(x)\equiv w(x^{\lambda_1})$ in $\Sigma_{\lambda_1}$ for some $\lambda_1\geq 0$, 
then $w$ is symmetric in the $x_1$ direction with respect to $T_{\lambda_1}$ and $\partial w/\partial x_1<0$ in $\{x_1>\lambda_1\}$. 
If this is not the case then
\begin{equation}
w(x_1,x')>w(-x_1,x')\ \mbox{for $x_1<0$},\quad 
\frac{\partial w}{\partial x_1}(x_1,x')<0 \ \mbox{for $x_1>0$}. 
 \label{eqn:pre-1-9}
\end{equation}
We again perform the procedure developed above for the negative $x_1$-direction to have either 
\[
w(x)=w(x^{\lambda_2})\ \mbox{in $\Sigma_{\lambda_2}$},\quad 
\frac{\partial w}{\partial x_1}>0\ \mbox{in $\{x_1<\lambda_2\}$}
\]
for some $\lambda_2\leq 0$, or 
\[
w(x_1,x')<w(-x_1,x')\ \mbox{for $x_1<0$},\quad 
\frac{\partial w}{\partial x_1}(x_1,x')>0 \ \mbox{for $x_1>0$}, 
\]
However, the latter and (\ref{eqn:pre-1-9}) do not simultaneously occur.  
Consequently, for some $\lambda\in{\bf R}$, $w$ is symmetric in the $x_1$ direction with respect to $T_{\lambda}$, 
$\partial w/\partial x_1<0$ in $\{x_1>\lambda\}$, and $\partial w/\partial x_1>0$ in $\{x_1<\lambda\}$. 

Since the problem is invariant with respect to rotation and translation, we can take any direction as the $x_1$-direction, 
and hence the lemma is established. \qed \\

When the compactness arguments are developed below, we often use the following lemma. 
\begin{lem}\label{lem:pre-2}
Given $R>0$, $C>0$, $A=A(x)\in C(\overline{B_R})$, let $v=v(x)$ be a solution of 
\[
\begin{cases}
-\Delta v=A(x)v_+^\gamma,\quad v\leq C & \mbox{in $B_R$} \\
v(x_0)=1 & \mbox{for some $x_0\in B_{R/2}$} 
\end{cases}
\]
Then there exists $C_{\const\insa}=C_{\suba}(n,\gamma,R,C,\|A\|_{L^\infty(B_R)})>0$ such that 
\[
v\geq -C_{\suba} \quad \mbox{in $B_{R/4}$}. 
\]
\end{lem}

{\it Proof.}\  
Let $v_1=v_1(x)$ and $v_2=v_2(x)$ be the solutions of 
\[
\begin{cases}
-\Delta v_1=A(x)v_+^\gamma & \mbox{on $B_R$} \\
v_1=0 & \mbox{on $\partial B_R$}, 
\end{cases}
\quad \mbox{and} \quad 
\begin{cases}
-\Delta v_2=0 & \mbox{on $B_R$} \\
v_2=v & \mbox{on $\partial B_R$}, 
\end{cases}
\]
respectively, i.e., $v=v_1+v_2$. 
Then the maximum principle and the representation formula imply 
\begin{equation}
0\leq v_1 \leq \|A\|_{L^\infty(B_R)}C^\gamma\sup_{x\in B_R}\|G_R(x,\cdot)\|_{L^1(B_R)}\equiv C_{\const\insz} \quad \mbox{in $B_R$}, 
 \label{eqn:pre-2-1}
\end{equation}
where $G_R=G_R(x,y)$ is the Green function associated to $-\Delta$ with the Dirichlet boundary condition in $B_R$. 
Moreover, it follows from the maximum principle that $v_2\leq C$ in $B_R$. 
We use 
\[
\max_{\overline{B_{R/2}}}v_2\geq v(x_0)-\max_{\overline{B_{R/2}}}v_1\geq 1-C_{\subz}
\]
and the Harnack inequality to the nonnegative harmonic function $C-v_2$, 
so that there exists $C_{\const\insy}=C_{\suby}(n)>0$ such that 
\[
\max_{\overline{B_{R/4}}}(C-v_2)\leq C_{\suby}\min_{\overline{B_{R/2}}}(C-v_2)
\leq C_{\suby}\{C-(1-C_\subz)\}, 
\]
or 
\begin{equation}
\min_{\overline{B_{R/4}}}v_2\geq C-C_{\suby}\{C-(1-C_\subz)\}. 
 \label{eqn:pre-2-2}
\end{equation}
The lemma follows from (\ref{eqn:pre-2-1}) and (\ref{eqn:pre-2-2}). \qed
\section{Proof of Theorem \ref{thm:class}}\label{sec:class}
We begin with 
\begin{lem}\label{lem:1-1}
There exist $C_0=C_0(n,\gamma)>0$ and $\delta_0=\delta_0(n,\gamma)>0$ such that 
\begin{equation}
\max_{\overline{B_{1/4}}} v \leq C_0 
\end{equation}
for any solution $v\in C^2(B_1)$ of 
\begin{equation}
\begin{cases}
-\Delta v = v_+^\gamma\quad\mbox{in $B_1$} \\ 
\int_{B_1} v_+^{\frac{n(\gamma-1)}{2}} < \delta_0 
\end{cases}
\end{equation} 
\end{lem}

{\it Proof.}\ 
We comply \cite{wy03}. 
Suppose that the assertion fails. 
Then there exists a solution sequence $v_k=v_k(x)$ of 
\begin{equation}
\begin{cases}
-\Delta v_k = (v_k)_+^\gamma & \mbox{in $B_1$} \\ 
\int_{B_1} (v_k)_+^{\frac{n(\gamma -1)}{2}} dx \leq \frac{1}{k} \\ 
\max_{\overline{B_{1/4}}} v_k \geq k.  
\end{cases}
 \label{eqn:1-1-1}
\end{equation}
For each $k$, we define $h_k \in C^2(B_1)$ and $y_k \in B_{1/2}$ by 
\[
h_k(y) = \left( \frac{1}{2} -r \right)^q v_k(y), \quad 
h_k(y_k) = \max_{\overline{B_{1/2}}} h_k(y),  
\]
where $q=\frac{2}{\gamma-1}$ and $r=|y|$. 
It holds that 
\begin{align}
h_k(y_k)&=\left( \frac{1}{2}-|y_k| \right)^q v_k(y_k) \geq \max_{\overline{B_{1/4}}} \left( \frac{1}{2} -r \right)^q v_k(y) \nonumber\\ 
&\geq \left( \frac{1}{4} \right)^q \max_{\overline{B_{1/4}}} v_k(y) \geq \left( \frac{1}{4} \right)^q k 
 \label{eqn:1-1-2}
\end{align}
for any $k$. 
Here we introduce 
\[
w_k(y) = \mu_k^q v_k(\mu_k y +y_k), 
\]
where 
\[
\sigma_k=\frac{1}{2}-|y_k|, \quad d_k^q=h_k(y_k)=\sigma_k^qv_k(y_k), \quad \mu_k = \sigma_k /d_k. 
\]
Since 
\[
\frac{1}{2}-|y| \geq \frac{1}{2} - (|y_k| + |y-y_k|) = \left( \frac{1}{2}-|y_k| \right) - |y-y_k| 
\geq \sigma_k - \frac{\sigma_k}{2} = \frac{\sigma_k}{2} 
\]
for any $y\in B_{\sigma_k /2}(y_k)$, it holds that 
\begin{equation}
d_k^q = h_k(y_k) \geq \left( \frac{1}{2} - |y| \right)^q v_k(y) \geq \left( \frac{\sigma_k}{2} \right)^q v_k(y) 
 \label{eqn:1-1-3} 
\end{equation}
for any $y\in B_{\sigma_k /2}(y_k)$. 

We use (\ref{eqn:1-1-1}) and (\ref{eqn:1-1-3}) to get 
\begin{equation}
\begin{cases}
-\Delta w_k = (w_k)_+^\gamma,\quad w_k \leq 2^q\quad  \mbox{in $B_{d_k/2}$} \\ 
\int_{B_{d_k /2}} (w_k)_+^{\frac{n(\gamma -1)}{2}} dx = \int_{B_{\sigma_k /2}(y_k)} (v_k)_+^{\frac{n(\gamma -1)}{2}} dx \leq \frac{1}{k} \\ 
w_k(0) = \mu_k^q v_k(y_k) = 1 
\end{cases} 
 \label{eqn:1-1-4} 
\end{equation} 
Note that $d_k\rightarrow+\infty$ by (\ref{eqn:1-1-2}). 
Hence the standard compactness argument, by Lemma \ref{lem:pre-2} and the elliptic regularity, admits $\tilde{w} \in C^2({\bf R}^n)$ 
such that $w_k \rightarrow \tilde{w}$ in $C_{loc}^{1+\alpha}({\bf R}^n)$ ($\alpha\in (0,1)$) and 
\[
\begin{cases}
-\Delta \tilde{w} = 0,\quad \tilde{w} \leq 2^q & \mbox{in ${\bf R}^n$} \\ 
\tilde{w}(0) = 1. &
\end{cases}
\]
Since $\tilde{w}=\tilde{w}(x)$ is harmonic and bounded above in ${\bf R}^n$, 
we have $\tilde{w} \equiv 1$ in ${\bf R}^n$ by the Liouville theorem. 
Therefore $w_k \rightarrow 1$ in $C_{loc}({\bf R}^n)$, which contradicts the second of (\ref{eqn:1-1-4}). \qed \\

At this stage, we can show 

\begin{lem}\label{lem:1-2}
Any classical solution of (\ref{eqn:higher-d'}) is bounded above. 
\end{lem}

{\it Proof.}\ 
Let $v=v(x)$ be a classical solution of (\ref{eqn:higher-d'}). 
Then there exists $R>0$ such that 
\[
\int_{{\bf R}^n \setminus B_R} v_+^{\frac{n(\gamma-1)}{2}} < \delta_0 
\]
by the constraint of (\ref{eqn:higher-d'}), where $\delta_0$ is as in Lemma \ref{lem:1-1}. 
Therefore, Lemma \ref{lem:1-1} gives  
\[
\sup_{{\bf R}^n \setminus B_{R+1}} v \leq C_0, 
\]
where $C_0$ is as in Lemma \ref{lem:1-1}. 
Noting that $\delta_0$ and $C_0$ are independent of $v=v(x)$, we obtain the lemma. \qed \\

By virtue of Lemma \ref{lem:1-2}, we can derive the representation formula of entire solutions via the Newton potential. 

\begin{lem}\label{lem:1-3}
For any nontrivial classical solution $v=v(x)$ of (\ref{eqn:higher-d'}), there exist $c_\gamma>0$ and $c_\gamma'>0$ such that 
\begin{align}
& v(x) = \frac{1}{(n-2)\omega_{n-1}} \int_{{\bf R}^n} |x-y|^{2-n} v_+^\gamma(y) dy - c_\gamma
 \label{eqn:lem:1-2-1} \\
& v(x) = -c_\gamma + c_\gamma' |x|^{2-n} + o(|x|^{2-n}) \quad\mbox{as $|x|\rightarrow+\infty$}. 
 \label{eqn:lem:1-2-2}
\end{align}
\end{lem}

{\it Proof.}\ 
Define $w=w(x)$ by 
\[
0 \leq w(x) = \frac{1}{(n-2)\omega_{n-1}} \int_{{\bf R}^n} |x-y|^{2-n} v_+^\gamma(y) dy. 
\]
At first, we shall show that $w=w(x)$ is well-defined and 
\begin{equation}
w(x)\rightarrow 0\quad \mbox{as $|x|\rightarrow+\infty$}.   
 \label{eqn:1-3-1}
\end{equation}
Lemma \ref{lem:1-2} and the constraint of (\ref{eqn:higher-d'}) assure 
\begin{equation}
v_+ \in L^s({\bf R}^n) \quad \mbox{for any $s \in \left[ \frac{n(\gamma-1)}{2}, \infty \right]$}. 
 \label{eqn:1-3-2}
\end{equation}
Given $R>0$, we introduce
\[
w_1(x) = \int_{|y-x|\geq R} |x-y|^{2-n}v_+^\gamma(y) dy, \quad w_2(x) = \int_{|y-x|< R} |x-y|^{2-n}v_+^\gamma(y) dy, 
\]
that is, $w=\frac{1}{(n-2)\omega_{n-1}}(w_1+w_2)$. 
Since $\gamma(n-1) \in [ n(\gamma-1), \infty)$, it holds by (\ref{eqn:1-3-2}) that 
\begin{align}
0 &\leq w_2(x) 
\leq \left( \int_{|z|<R} |z|^{1-n}dz \right)^{\frac{n-2}{n-1}} \left( \int_{|z|<R} v_+^{\gamma(n-1)}(x-z)dz \right)^{\frac{1}{n-1}} \nonumber\\ 
&\leq C_{\const}(n,R) \| v_+ \|_{L^{\gamma(n-1)}(B(x,R))}^\gamma \rightarrow 0 \quad \mbox{as $|x| \rightarrow +\infty$}. 
 \label{eqn:1-3-3}
\end{align}
In addition, $w_1$ is estimated by 
\begin{align}
0 \leq w_1(x) &\leq 
 \begin{cases}
 &R^{2-n} \int_{|z|\geq R} v_+^\gamma (x-z) dz \quad \mbox{if $\gamma \in \left(1, \frac{n}{n-2} \right]$} \\ 
 &\left( \int_{|z|\geq R} |z|^{-n \left(1+\frac{2}{(n-2)\gamma -n}\right)} dz \right)^{\frac{(n-2)\gamma -n}{n(\gamma -1)}} \\ 
 &\quad\times \left( v_+^{\frac{n(\gamma -1)}{2}} dz \right)^{\frac{2\gamma}{n(\gamma -1)}} 
 \quad \mbox{if $\gamma \in \left(\frac{n}{n-2}, \frac{n+2}{n-2} \right)$}
 \end{cases} \nonumber\\
&\leq 
 \begin{cases}
 &R^{2-n} \| v_+ \|_\gamma^\gamma \quad \mbox{if $\gamma \in \left(1, \frac{n}{n-2} \right]$} \\ 
 &R^{-\frac{1}{\gamma -1}} C_{\const}(n,\gamma) \| v_+ \|_{\frac{n(\gamma -1)}{2}}^\gamma  
 \quad \mbox{if $\gamma \in \left(\frac{n}{n-2}, \frac{n+2}{n-2} \right)$}
 \end{cases}
 \label{eqn:1-3-4}
\end{align} 
From (\ref{eqn:1-3-2})-(\ref{eqn:1-3-4}) and the property that $\gamma\in \left[ \frac{n(\gamma-1)}{2}, \infty \right)$ for $\gamma \in \left(1, \frac{n}{n-2} \right]$, 
we see that $w$ is well-defined and 
\[
0 \leq \limsup_{|x| \rightarrow +\infty} w(x) \leq 
 \begin{cases}
 C_{\const}(n,\gamma) R^{2-n} & \mbox{if $\gamma \in \left(1, \frac{n}{n-2} \right]$} \\ 
 C_{\const}(n,\gamma) R^{-\frac{1}{\gamma -1}} & \mbox{if $\gamma \in \left(\frac{n}{n-2}, \frac{n+2}{n-2} \right)$,}
 \end{cases}
\]
which implies (\ref{eqn:1-3-1}) since $R >0$ is arbitrary. 

Next we have 
\[
-\Delta(v-w)=0 \quad \mbox{in ${\bf R}^n$}, \quad \sup_{{\bf R}^n} (v-w) < +\infty  
\]
by Lemma \ref{lem:1-2} and (\ref{eqn:1-3-1}). 
Then the Liouville theorem guarantees that there exists $c_1 \in {\bf R}$ such that $v-w=c_1$. 
Note that $c_1<0$, since if this is not the case then $-\Delta v =v^{\gamma},\ v \geq 0$ in ${\bf R}^n$, 
which is impossible from $\gamma\in\left(1,\frac{n+2}{n-2}\right)$ and the result of \cite{gs81-sub}, 
recall that $v$ is now a nontrivial classical solution of (\ref{eqn:higher-d'}). 
Hence we obtain (\ref{eqn:lem:1-2-1}) for $c_\gamma=-c_1>0$. 

Finally, (\ref{eqn:lem:1-2-2}) holds for $c_\gamma' =  \frac{1}{(n-2)\omega_{n-1}} \int_{{\bf R}^n} v_+^\gamma dx$ since  
\begin{align*}
|x|^{n-2} (v(x)+c_\gamma) &=|x|^{n-2}w(x)=\frac{1}{(n-2)\omega_{n-1}} \int_{{\bf R}^n} \frac{|x|^{n-2}}{|x-y|^{n-2}} v_+^\gamma(y) dy \\ 
&\rightarrow \frac{1}{(n-2)\omega_{n-1}} \int_{{\bf R}^n} v_+^\gamma dx \quad \mbox{as $|x| \rightarrow +\infty$}, 
\end{align*}
by (\ref{eqn:1-3-1})-(\ref{eqn:1-3-2}) and the dominated convergence theorem. \qed \\

Now we prove Theorem \ref{thm:class}. \\

\underline{{\it Proof of Theorem \ref{thm:class}}}\ \ Let $v=v(x)$ be a nontrivial classical solution of (\ref{eqn:higher-d'}), and consider 
\[
0\leq w(x)=\frac{1}{(n-2)\omega_{n-1}} \int_{{\bf R}^n} |x-y|^{2-n} v_+^\gamma(y) dy. 
\]
From Lemma \ref{lem:1-3} and its proof, we see that $w=v+c_\gamma$ for some $c_\gamma>0$, and that  
\[
\begin{cases}
-\Delta w=(w-c_\gamma)_+^\gamma,\ w>0 & \mbox{in ${\bf R}^n$} \\
w(x)\rightarrow 0 & \mbox{as $|x|\rightarrow+\infty$}. 
\end{cases}
\]
We apply Lemma \ref{lem:pre-1} and conclude that $w=w(x)$ is radially symmetric about some $x_0\in{\bf R}^n$ and satisfies $\partial u/\partial r<0$ for $r=|x-x_0|>0$, and so is and does $v=v(x)$. 

It is left to prove that $v=v(x)$ is represented by (\ref{eqn:thm:class-1}). 
The remainders of the statement of the theorem directly follow from (\ref{eqn:thm:class-1}). 
Let $\phi=\phi(r)$ be the unique classical solution to (\ref{eqn:phi-ode}). 
In view of the radial symmetry of $v=v(x)$, we put $v(x)=\psi(r)$, where $r=|x-x_0|$. 
Noting the scaling invariance of the problem, we see that $\psi_0(r)=\mu_0^q\psi(\mu_0 r)$ satisfies (\ref{eqn:phi-ode}), where $\mu_0=\psi(0)^{-1/q}$ and $q=\frac{2}{\gamma-1}$. 
Hence the uniqueness of the problem (\ref{eqn:phi-ode}) assures that $\psi_0(r)=\phi(r)$ for $r\geq 0$, in particular, 
\[
\psi_0(r)=\frac{\lambda_\gamma^\ast}{\omega_{n-1}(n-2)}\left(\frac{1}{r^{n-2}}-\frac{1}{(r_\gamma^\ast)^{n-2}}\right) \quad \mbox{for $r\geq r_\gamma^\ast$}, 
\]
where $r_\gamma^\ast>0$ is the first zero point of $\phi=\phi(r)$.
The proof is complete. \qed
\section{Proof of Theorem \ref{thm:supinf}}\label{sec:supinf}
Theorem \ref{thm:supinf} is a direct consequence of the following lemma. 

\begin{lem}\label{lem:2-1}
Assume the assumptions of Theorem \ref{thm:supinf} for $\Omega=B_1$. 
Then, given $T>0$, there exist $C_{\const\insa}=C_{\suba}(n,\gamma,A)>0$ and $C_{\const\insb}=C_{\subb}(n,\gamma,A,T)>0$ such that 
\begin{equation}
v(0) + C_{\suba} \inf_{B_1} v \leq C_{\subb} 
 \label{eqn:lem:2-1-1}
\end{equation}
for any solution $v=v(x)$ of 
\begin{equation}
\begin{cases} 
-\Delta v = A(x)v_+^\gamma \quad \mbox{in $B_1$} \\ 
\int_{B_1} v_+^{\frac{n(\gamma -1)}{2}} dx \leq T. 
\end{cases}
 \label{eqn:lem:2-1-2}
\end{equation}
\end{lem}

{\it Proof.}\ 
Assume that the statement is false.  
Then, given $\hat{C} >0$, there exists a solution sequence $v_k=v_k(x)$ of 
\begin{equation}
\begin{cases}
-\Delta v_k = A(x)(v_k)_+^\gamma \quad \mbox{in $B_1$} \\ 
\int_{B_1} (v_k)_+^{\frac{n(\gamma -1)}{2}} dx \leq T \\ 
v_k(0) + \hat{C} \inf_{B_1} v_k \geq k. 
\end{cases}
 \label{eqn:2-1-1}
\end{equation}
It is obvious that 
\begin{equation*}
v_k(0) \geq \frac{k}{1+\hat{C}} \rightarrow +\infty \quad \mbox{as $k \rightarrow \infty$}. 
\end{equation*} 
As in the proof of Lemma \ref{lem:1-1}, we introduce 
\begin{align*}
&w_k(y) = \mu_k^q v_k(\mu_k y +y_k),\quad h_k(y_k)=\max_{\overline{B_{1/2}}} h_k(y),\quad h_k(y) = \left( \frac{1}{2} -r \right)^q v_k(y), \\
&\sigma_k=\frac{1}{2}-|y_k|, \quad d_k^q=h_k(y_k)=\sigma_k^qv_k(y_k), \quad \mu_k = \sigma_k /d_k, 
\end{align*}
and find
\begin{align}
&d_k \geq \frac{1}{2}v_k(0)^{1/q} \rightarrow +\infty, 
 \label{eqn:2-1-3} \\
&w_k \leq 2^q \quad \mbox{in $B_{d_k /2}(y_k)$}. 
 \label{eqn:2-1-4}
\end{align} 
From the scaling invariance, (\ref{eqn:2-1-1}) and (\ref{eqn:2-1-4}), it follows that 
\begin{equation*}
\begin{cases}
-\Delta w_k = A_k'(y)(w_k)_+^\gamma,\quad w_k \leq 2^q\quad \mbox{in $B_{d_k/2}$} \\ 
\int_{B_{d_k /2}} (w_k)_+^{\frac{n(\gamma -1)}{2}} dx = \int_{B_{\sigma_k /2}(y_k)} (v_k)_+^{\frac{n(\gamma -1)}{2}} dx \leq T \\ 
w_k(0) = 1 
\end{cases} 
\end{equation*} 
where $A_k'(y)=A(\mu_k y+y_k)$. 
We may assume that $y_k\rightarrow y_0$ for some $y_0\in \overline{B_{1/2}}$. 
Then there exists $\tilde{w}\in C^2({\bf R}^n)$ such that $w_k\rightarrow\tilde{w}$ in $C_{loc}^{1+\alpha}({\bf R}^n)$ ($\alpha\in (0,1)$) and 
\begin{equation}
\begin{cases}
-\Delta \tilde{w} = A_0\tilde{w}_+^\gamma,\quad \tilde{w} \leq 2^q\quad \mbox{in ${\bf R}^n$} \\ 
\int_{{\bf R}^n} \tilde{w}_+^{\frac{n(\gamma -1)}{2}} dx \leq T \\ 
\tilde{w}(0) = 1 
\end{cases}
 \label{eqn:2-1-5}
\end{equation}
by (\ref{eqn:2-1-3}) and Lemma \ref{lem:pre-2}, where $A_0=A(y_0)$. 

Here we note that $A_0=A(y_0)>0$. 
Actually, if this is not the case (i.e., $A_0=0$) then the Liouville theorem leads to a contradiction by the integrable condition of (\ref{eqn:2-1-5}).  

Now we put 
\[
\tilde{z}(x)=A_0^{\frac{1}{\gamma-1}}\tilde{w}(x), \quad z_k(x)=A_0^{\frac{1}{\gamma-1}}w_k(x), 
\]
and see that $z_k\rightarrow\tilde{z}$ in $C_{loc}^{1+\alpha}({\bf R}^n)$ ($\alpha\in (0,1)$) and 
\begin{equation}
\begin{cases}
-\Delta \tilde{z} = \tilde{z}_+^\gamma,\quad \tilde{z} \leq 2^q A_0^{\frac{1}{\gamma-1}}\quad \mbox{in ${\bf R}^n$} \\ 
\int_{{\bf R}^n} \tilde{z}_+^{\frac{n(\gamma -1)}{2}} dx \leq A_0^{n/2}T \\ 
\tilde{z}(0) = A_0^{\frac{1}{\gamma-1}}.
\end{cases}
 \label{eqn:2-1-6}
\end{equation}
By virtue of Theorem \ref{thm:class}, $\tilde{z}=\tilde{z}(x)$ is represented by (\ref{eqn:thm:class-1}) for some $x_0\in{\bf R}^n$ and $\mu\in[\sqrt{A_0},2\sqrt{A_0}]$, in particular, 
\[
\tilde{z}(x)\rightarrow -\sigma_0 \quad \mbox{as $|x|\rightarrow +\infty$}
\]
for some $\sigma_0>0$. 
Therefore, there exist $C_{\const\insc}=C_{\subc}(n,\gamma,A)>0$ and $R=R(n,\gamma,A)>0$, independent of $\hat{C}$, such that 
\[
\tilde{z}(0)+C_{\subc}\inf_{\partial B_R}\tilde{z}<0, 
\]
or 
\begin{equation}
w_k(0)+C_{\subc}\inf_{\partial B_R}w_k<0. 
 \label{eqn:2-1-7}
\end{equation} 
for $k \gg 1$. 
Noting that $v_k$ is super-harmonic, we obtain
\[
v_k(0)+C_{\subc}\inf_{B_1}v_k \leq v_k(y_k)+C_{\subc}\inf_{\partial B_{\mu_k R}(y_k)}v_k=\mu_k^{-q}\left(w_k(0)+C_{\subc}\inf_{\partial B_R} w_k \right)<0  
\]
for $k \gg 1$ by (\ref{eqn:2-1-7}), 
which contradicts (\ref{eqn:2-1-1}) when $\hat{C}=C_{\subc}$ since $\hat{C}>0$ is arbitrary and $C_{\subc}$ is independent of $\hat{C}$. \qed \\

\underline{{\it Proof of Theorem \ref{thm:supinf}}}\ \ There exist $\mu_0=\mu_0(K) >0$ and $x_0 \in K$ such that 
\[
\bigcup_{x \in K} B_{\mu_0}(x) \subset \Omega, \quad v(x_0)=\sup_{K} v.   
\]
We introduce $w(x)=\mu_0^q v(\mu_0 x+x_0)$ for $x \in B_1$ and $q=\frac{2}{\gamma -1}$. 
Then $w=w(x)$ satisfies 
\[
\begin{cases} 
-\Delta w=A(\mu_0 x+x_0)w_+^\gamma \quad \mbox{in $B_1$} \\ 
\int_{B_1} w_+^{\frac{n(\gamma -1)}{2}} dx \leq T 
\end{cases}
\]
and 
\begin{equation}
v(x_0)+C \inf_{\Omega}v \leq v(x_0)+C \inf_{B_{\mu_0}(x)} v=\mu_0^{-q}(w(0)+C \inf_{B_1} w) 
 \label{eqn:2-1-8}
\end{equation}
for any $C>0$. 
On the other hand, Lemma \ref{lem:2-1} admits $C_{\const\insd}=C_{\subd}(n,\gamma,A)>0$ and $C_{\const\inse}=C_{\sube}(n,\gamma,A,T)>0$ such that 
\begin{equation}
w(0)+C_{\subd}\inf_{B_1}w\leq C_{\sube}. 
 \label{eqn:2-1-9}
\end{equation}
The theorem follows from (\ref{eqn:2-1-8}) and (\ref{eqn:2-1-9}). \qed
\section{Proof of Theorem \ref{thm:quantization}}\label{sec:bm-ls-type} 
In this section, we use the notation
\[
A_{\max}=\max_{\overline{\Omega}}A. 
\]
We start with the $\varepsilon$-regularity. 

\begin{lem}\label{lem:epsilon-reg}
Assume the assumptions of Theorem \ref{thm:quantization} for $\Omega=B_R$, $R>0$. 
Then, given $\varepsilon\in(0,A_{\max}^{-n/2}\lambda_\gamma^\ast)$, there exists $C_{\varepsilon,R}=C_{\varepsilon,R}(n,\gamma,A_k,A,\varepsilon,R)>0$ such that 
\[
\max_{\overline{B_{R/4}}}v_k\leq C_{\varepsilon,R}
\]
for any solution sequence $v_k=v_k(x)$ of 
\[
\begin{cases} 
-\Delta v_k = A_k(x)(v_k)_+^\gamma \quad \mbox{in $B_R$} \\ 
\int_{B_R} (v_k)_+^{\frac{n(\gamma -1)}{2}} dx \leq \varepsilon. 
\end{cases}
\]
\end{lem}

{\it Proof.}\ 
We have only to show the lemma for the case $R=1$ thanks to the scaling invariance. 
Fix $\varepsilon\in(0,A_{\max}^{-n/2}\lambda_\gamma^\ast)$ and assume that the assertion fails. 
Then, similarly to the proof of Lemma \ref{lem:2-1}, we obtain $\tilde{w}=\tilde{w}(x)\in C^2({\bf R}^n)$ such that 
\[
\begin{cases}
-\Delta \tilde{w} = A_0\tilde{w}_+^\gamma,\quad \tilde{w} \leq 2^q\quad \mbox{in ${\bf R}^n$} \\ 
\int_{{\bf R}^n} \tilde{w}_+^{\frac{n(\gamma -1)}{2}} dx \leq \varepsilon \\ 
\tilde{w}(0) = 1 
\end{cases}
\]
for some $0<A_0\leq A_{\max}$. 
Then $\tilde{z}(x)=A_0^{1/(\gamma-1)}\tilde{w}(x)$ satisfies  
\[
\begin{cases}
-\Delta \tilde{z}=\tilde{z}_+^\gamma,\quad \tilde{z} \leq 2^q A_0^{\frac{1}{\gamma-1}}\quad \mbox{in ${\bf R}^n$} \\ 
\int_{{\bf R}^n} \tilde{z}_+^{\frac{n(\gamma -1)}{2}} dx \leq A_0^{n/2}\varepsilon \leq A_{\max}^{n/2}\varepsilon<\lambda_\gamma^\ast \\ 
\tilde{z}(0)=A_0^{\frac{1}{\gamma-1}}, 
\end{cases}
\]
however, there is no such a solution by Theorem \ref{thm:class}. \qed \\

The proof of Theorem \ref{thm:quantization} is reduced to showing the following two propositions. 

\begin{pro}\label{pro:bm}
Assume that $\gamma \in \left(1, \frac{n+2}{n-2} \right)$, $n \geq 3$, $\Omega$ is an open set, 
$0\leq A_k\in C(\overline{\Omega})$, 
and there exists $0\leq A\in C(\overline{\Omega})$ such that $A_k\rightarrow A$ in $C(\overline{\Omega})$. 
Given $T>0$, let $\{v_k\}$ be a solution sequence of 
\begin{equation}
\begin{cases}
-\Delta v_k = A_k(x)(v_k)_+^\gamma \quad \mbox{in $\Omega$} \\ 
\int_{\Omega} (v_k)_+^{\frac{n(\gamma -1)}{2}} dx \leq T. 
\end{cases}
 \label{eqn:pro:bm-0}
\end{equation}
Then, passing to a subsequence, we have the following alternatives: 

{\bf (i)} $\{ v_k \}$ is locally uniformly bounded in $\Omega$. 

{\bf (ii)} $v_k \rightarrow -\infty$ locally uniformly in $\Omega$. 

{\bf (iii)} There exists a finite set ${\cal S} = \{x_i\}_{i=1}^l$ such that 
$v_k \rightarrow -\infty$ locally uniformly in $\Omega\setminus{\cal S}$, and that 
\begin{equation}
(v_k)_+^{\frac{n(\gamma -1)}{2}}dx \overset{\ast}{\rightharpoonup} \sum_{i=1}^{l} m(x_i) \delta_{x_i}(dx) \quad \mbox{in ${\cal M}(\Omega)$}
 \label{eqn:pro:bm-1}
\end{equation}
with $m(x_i)\geq A_{\max}^{-n/2}\lambda_\gamma^\ast$ for $i=1,\cdots,l$. 
\end{pro}

\begin{pro}\label{pro:ls}
It holds that $m(x_i)\in A(x_i)^{-n/2}\lambda_\gamma^\ast{\bf N}$ for each $i=1,\cdots,l$ in (\ref{eqn:pro:bm-1}). 
\end{pro}

We first give the proof of Proposition \ref{pro:bm}. \\

\underline{{\it Proof of Proposition \ref{pro:bm}}}\ \ Since $\{ (v_k)_+^{\frac{n(\gamma-1)}{2}} \}$ is bounded in $L^1(\Omega)$, 
we have a non-negative measure $\mu$ such that 
\[
(v_k)_+^{\frac{n(\gamma-1)}{2}} dx \overset{\ast}{\rightharpoonup} \mu \quad \mbox{in ${\cal M}(\Omega)$}, 
\]
passing to a subsequence. 
We put 
\begin{align}
&{\cal S} = \{ x \in\Omega \ | \ \mbox{there is $\{ x_k \} \subset \Omega$ such that $x_k \rightarrow x$ and $v_k(x_k) \rightarrow +\infty$} \} 
 \label{eqn:ls-S} \\ 
&\Sigma = \{ x \in\Omega \ | \ \mu(\{x\}) \geq A_{\max}^{-n/2}\lambda_\gamma^\ast \}. 
 \label{eqn:ls-Sigma} 
\end{align}

First we shall show that ${\cal S}=\Sigma$. 
If $x_0\not\in\Sigma$ then there exists $r_0>0$ such that $\mu(B_{r_0}(x_0))<A_{\max}^{-n/2}\lambda_\gamma^\ast$, 
and hence
\[
\int_{B_{r_0}(x_0)} (v_k)_+^{\frac{n(\gamma-1)}{2}} dx \leq \varepsilon_0
\]
for $k\gg 1$ and for some $\varepsilon_0\in(0,A_{\max}^{-n/2}\lambda_\gamma^\ast)$. 
Thus Lemma \ref{lem:epsilon-reg} assures that $x_0\not\in{\cal S}$. 
In turn, if $x_0\not\in{\cal S}$ then there exists $r_1>0$ such that $\sup_k\|(v_k)_+\|_{L^\infty(B_{r_1}(x_0))}<+\infty$, and hence 
\[
\lim_{r \downarrow 0} \limsup_{k\rightarrow\infty} \int_{B(x_0,r_0)} (v_k)_+^{\frac{n(\gamma-1)}{2}} dx = 0. 
\]
This means that $x_0\not\in\Sigma$, and therefore ${\cal S}=\Sigma$. 

Next we shall show that either (i) or (ii) occurs if ${\cal S} = \emptyset$. 
Fix an open set $\omega$ such that $\overline{\omega} \subset \Omega$ and $\overline{\omega}$ is compact. 
There exists $C_{\const\insa}>0$ such that 
\[
\sup_k \| (v_k)_+ \|_{L^\infty(\omega)} \leq C_{\suba}. 
\]
Let $v_{1,k}$ be a solution of 
\[
\begin{cases}
-\Delta v_{1,k} = (v_k)_+^\gamma & \mbox{in $\omega$} \\ 
v_{1,k} = 0 & \mbox{on $\partial\omega$.}
\end{cases}
\] 
It holds that $v_{1,k} \geq 0$ in $\omega$ by the maximum principle, 
and that $\{ v_{1,k} \}$ is uniformly bounded in $\omega$ by the elliptic regularity. 
Hence $v_{2,k} = v_k - v_{1,k}$ is harmonic and bounded above in $\omega$,   
and then the Harnack principle guarantees that $\{v_{2,k}\}$ is uniformly bounded in $\omega$, 
otherwise $v_{2,k} \rightarrow -\infty$ in $\omega$. 
These alternatives hold for $v_k$ since $v_{1,k}$ is uniformly bounded in $\omega$. 
Since $\omega$ is arbitrary, either (i) or (ii) occurs if ${\cal S} = \emptyset$. 

Finally, we shall show that ${\cal S} \neq \emptyset$ implies (iii). 
The proof of this part is different from \cite{wy03}. 
We adopt the blowup analysis here. 
We may put ${\cal S}=\{x_i\}_{i=1}^l$ by ${\cal S}=\Sigma$ and $\mu(\Omega)\leq T$. 
Similarly to the argument above, we see that either (I) or (II) below holds: 

(I) $\{v_k\}$ is locally uniformly bounded in $\Omega\setminus{\cal S}$. 

(II) $v_k \rightarrow -\infty$ locally uniformly in $\Omega\setminus{\cal S}$. 

We claim that (I) cannot occur. 
Assume that (I) occur, and fix $x_0\in{\cal S}$. 
Then there exist $r_2>0$ and $C_{\const\insb}>0$ such that 
$B_{r_2}(x_0) \cap {\cal S} = \{x_0\}$ and 
\begin{equation}
v_k \geq -C_{\subb} \quad \mbox{on $\partial B_{r_2}(x_0)$.}
 \label{eqn:bm-1}
\end{equation}
It follows from the definition of ${\cal S}$ and (I) that there exists a maximizer $x_k$ of $v_k$ in $B_{r_2}(x_0)$. 
It is obvious that $x_k\rightarrow x_0$ by (I).  
Function $z_k=z_k(x)$ defined by 
\[
z_k(x)=A_0^{\frac{1}{\gamma-1}}\mu_k^q v_k(\mu_k x + x_k),\quad \mu_k^{-q}=v_k(x_k),\quad q=\frac{2}{\gamma-1},\quad A_0=A(x_0)>0 
\]
satisfies 
\[
\begin{cases}
-\Delta z_k = \frac{A_k''(x)}{A_0}(z_k)_+^\gamma & \mbox{in $B_{\frac{r_0}{2\mu_k}}$} \\ 
\int_{B_{\frac{r_0}{2\mu_k}}} (z_k)_+^{\frac{n(\gamma -1)}{2}} dx \leq A_0^{n/2}T & \\
-A_0^{\frac{1}{\gamma-1}}\mu_k^{q}C_{\subb} \leq z_k \leq z_k(0)=A_0^{\frac{1}{\gamma-1}} & \mbox{in $B_{\frac{r_0}{2\mu_k}}$} 
\end{cases}
\]
by the maximum principle, where $A_k''(x)=A_k(\mu_k x+x_k)$.  
Hence the compactness argument admits $z=z(x)\in C^2({\bf R}^n)$ such that $z_k \rightarrow z$ in $C_{loc}^{1+\alpha}({\bf R}^n)$ ($\alpha\in (0,1)$) and 
\[
\begin{cases}
-\Delta z = z_+^\gamma,\quad 0 \leq z \leq z(0)=A_0^{\frac{1}{\gamma-1}}\quad \mbox{in ${\bf R}^n$} \\ 
\int_{{\bf R}^n} z_+^{\frac{n(\gamma-1)}{2}} dx \leq A_0^{n/2}T, \\ 
\end{cases}
\]
which is impossible by Theorem \ref{thm:class}. 

We have shown that if ${\cal S} \neq \emptyset$ then $v_k \rightarrow -\infty$ locally uniformly in $\Omega\setminus{\cal S}$. 
Therefore, $(v_k)_+^{\frac{n(\gamma-1)}{2}} \rightarrow 0$ in $L_{loc}^1(\Omega\setminus{\cal S})$, and hence 
\[
\mu(dx) = \sum_{x_0 \in {\cal S}} m(x_0) \delta_{x_0}(dx) 
\]
with $m(x_0)\geq A_{\max}^{-n/2}\lambda_\gamma^\ast$ for any $x_0\in {\cal S}$. 
The proof is complete. \qed \\

Here we prepare the key estimate to prove Proposition \ref{pro:ls}. 
The proof is done similarly to \cite{wy03}. 

\begin{lem}\label{lem:3-0}
Assume the assumption of Theorem \ref{thm:quantization} for $\Omega=B_R$, $R>0$. 
Let $v_k=v_k(x)$ be a solution sequence of 
\begin{align}
&-\Delta v_k=A_k(x)(v_k)_+^\gamma \quad \mbox{in $B_R$} 
 \label{eqn:sharp}\\ 
&\int_{B_R} (v_k)_+^{\frac{n(\gamma-1)}{2}}dx\leq T 
 \label{eqn:sharp'}\\
&v_k(x)|x|^q \leq C_{\const\inse} \quad \mbox{for $x\in B_R\setminus B_{R_0/2}$}, 
 \label{eqn:rate}
\end{align}
where $C_{\sube}>0$ and $R_0\in (0,R/2)$. 
Then, there exist $C_i=C_i(n,\gamma,T,A_k,A,C_{\sube})>0$ ($i=\const\insg,\const\insh$) such that 
\begin{equation}
\sup_{\partial B_r}v_k \leq -C_{\subg}v_k(0) +r^{-q}C_{\subh}
 \label{eqn:lem:3-2-2}
\end{equation}
for any $r\in[R_0,R/2)$, where $C_i$ ($i=\subg,\subh$) are independent of $v_k$, $R$, $R_0$ and $r$. 
\end{lem}

The proof of Lemma \ref{lem:3-0} is stated later on. 
For the purpose, we prepare the two estimates below. 

\begin{lem}\label{lem:3-1}
Assume the assumption of Theorem \ref{thm:quantization} for $\Omega=B_R$, $R>0$, and let $v_k=v_k(x)$ satisfy (\ref{eqn:sharp})-(\ref{eqn:sharp'}). 
Then, there exist $C_{\const\insc}=C_{\subc}(n,\gamma,A_k,A)>0$ and $C_{\const\insd}=C_{\subd}(n,\gamma,A_k,A,T)$ such that 
\[
v_k(0)+C_{\subc}\inf_{\partial B_r}v_k \leq r^{-q}C_{\subd} 
\]
for any $r\in (0,R)$, where $C_{\subc}$ and $C_{\subd}$ are independent of $R$ and $r$. 
\end{lem}

{\it Proof.}\ 
We put $v_k^{(r)}(x) = r^q v_k(rx)$ and $A_k^{(r)}(x)=A_k(rx)$ for $r\in (0,R)$ and $q=\frac{2}{\gamma-1}$, so that 
\begin{equation}
\begin{cases}
-\Delta v_k^{(r)} = A_k^{(r)}(x)(v_k^{(r)})_+^\gamma \quad \mbox{in $B_1$} \\ 
\int_{B_1} (v_k^{(r)})_+^{\frac{n(\gamma-1)}{2}} dx \leq T. 
\end{cases}
 \label{eqn:3-1-1}
\end{equation}
The argument developed in the proof of Lemma \ref{lem:2-1} still works for (\ref{eqn:3-1-1}), 
and therefore we obtain $C_{\const\insi}=C_{\subi}(n,\gamma,A_k,A)>0$ and $C_{\const\insj}=C_{\subj}(n,\gamma,A_k,A,T)>0$ such that 
\[
v_k^{(r)}(0)+C_{\subi} \inf_{B_1} v_k^{(r)}=v_k^{(r)}(0)+C_{\subi} \inf_{\partial B_1} v_k^{(r)}\leq C_{\subj}, 
\]
which yields the desired estimate. \qed \\

\begin{lem}\label{lem:3-2}
Assume the assumption of Theorem \ref{thm:quantization} for $\Omega=B_R$, $R>0$, and let $v_k=v_k(x)$ satisfy (\ref{eqn:sharp}) and (\ref{eqn:rate}). 
Then, there exist $C_{\const\insf}=C_{\subf}(n,\gamma,A_k,A,C_{\sube})>0$ and $\beta=\beta(n)\in(0,1)$ such that 
\[
\sup_{\partial B_r}v_k \leq \beta\inf_{\partial B_r}v_k + r^{-q} C_{\subf} 
\]
for any $r\in[R_0/2,R)$, where $C_{\subf}$ and $\beta$ are independent of $v_k$, $R$, $R_0$ and $r$. 
\end{lem}

{\it Proof.}\ 
Given $r\in [R_0, R/2)$, we put $v_k^{(r)}(x) = r^q v_k(rx)$ and $A_k^{(r)}(x)=A_k(rx)$, where $q=\frac{2}{\gamma-1}$. 
Then $v_k^{(r)}=v_k^{(r)}(x)$ satisfies 
\[
\begin{cases}
-\Delta v_k^{(r)}=A_k^{(r)}(x)(v_k^{(r)})_+^\gamma & \mbox{in $B_2 \setminus B_{1/2}$} \\ 
v_k^{(r)}\leq 2^q C_{\sube} & \mbox{on $\overline{B_2 \setminus B_{1/2}}$}. 
\end{cases}
\]
Let $w_k^{(r)}=w_k^{(r)}(x)$ be the solution of 
\[
\begin{cases}
-\Delta w_k^{(r)}=A_k^{(r)}(x)(v_k^{(r)})_+^\gamma & \mbox{in $B_2 \setminus B_{1/2}$} \\ 
w_k^{(r)}=0 & \mbox{on $\partial(B_2 \setminus B_{1/2})$}. 
\end{cases}
\]
Then there exists $C_{\const\insk}=C_{\subk}(n,\gamma,\sup_k\|A_k\|_{L^\infty(B_R)},C_{\sube})>0$ such that 
\begin{equation}
0 \leq w_k^{(r)} \leq C_{\subk} \quad \mbox{on $\overline{B_2 \setminus B_{1/2}}$}. 
 \label{eqn:3-2-1}
\end{equation}
by the maximum principle and the elliptic regularity. 

Since $\xi_k^{(r)}(x)=v_k^{(r)}(x)-w_k^{(r)}(x)$ is harmonic in $B_2 \setminus B_{1/2}$, 
and since $\xi_k^{(r)} \leq 2^q C_{\sube}$ on $\overline{B_2 \setminus B_{1/2}}$, 
$2^qC_{\sube}-\xi_k^{(r)}$ is nonnegative, bounded above and harmonic in $B_2 \setminus B_{1/2}$. 
Hence the Harnack inequality admits $\beta=\beta(n)\in(0,1)$, independent of $v_k$, $R$, $R_0$ and $r$, such that 
\begin{equation}
\beta \sup_{\partial B_1} (2^qC_{\sube}-\xi_k^{(r)}) \leq \inf_{\partial B_1} (2^qC_{\sube}-\xi_k^{(r)}).  
 \label{eqn:3-2-2}
\end{equation}
Inequalities (\ref{eqn:3-2-1}) and (\ref{eqn:3-2-2}) imply
\begin{align*}
&\sup_{\partial B_1} v_k^{(r)}=\sup_{\partial B_1} \{ -(2^qC_{\sube}-\xi_k^{(r)}) + 2^qC_{\sube} + w_k^{(r)} \} \\ 
&\leq -\inf_{B_1} (2^qC_{\sube}-\xi_k^{(r)}) + 2^qC_{\sube} + \sup_{\partial B_1} w_k^{(r)} \\ 
&\leq -\beta \sup_{\partial B_1} (2^qC_{\sube}-\xi_k^{(r)}) + 2^qC_{\sube} + C_{\subk}
=\beta \inf_{\partial B_1} (\xi_k^{(r)}-2^qC_{\sube}) + 2^qC_{\sube} + C_{\subk} \\ 
&\leq \beta \inf_{\partial B_1} v_k^{(r)} + 2^q(1-\beta)C_{\sube} + C_{\subk}, 
\end{align*}
and thus the lemma is shown. \qed \\

\underline{{\it Proof of Lemma \ref{lem:3-0}}} The lemma follows from Lemmas \ref{lem:3-1} and \ref{lem:3-2}. \qed \\ 

We readily see that the proof of Proposition \ref{pro:ls} is reduced to showing Lemmas \ref{lem:3y}-\ref{lem:3z} below. 


\begin{lem}\label{lem:3y}
Assume the assumptions of Theorem \ref{thm:quantization} for $\Omega=B_R$, $R>0$. 
Let $v_k=v_k(x)$ satisfy 
\begin{align}
&-\Delta v_k = A_k(x)(v_k)_+^\gamma \quad \mbox{in $B_R$} 
 \label{eqn:lem:3x-1} \\ 
&\max_{\overline{B_R}} v_k \rightarrow +\infty \quad \mbox{and} \quad \max_{\overline{B_R}\setminus B_r} v_k \rightarrow -\infty 
 \label{eqn:lem:3x-2} \\ 
&\int_{B_R} (v_k)_+^{\frac{n(\gamma-1)}{2}} dx \leq T 
 \label{eqn:lem:3y-0}
\end{align}
for any $k$ and $r\in(0,R)$, and for some $T>0$. 
Then, passing to a subsequence, we have $\{ x_k^{(j)} \}_{j=0}^{m-1} \subset B_R$, $\{ l_k^{(j)} \}_{j=0}^{m-1} \subset (0,+\infty)$ 
and $m \in {\bf N}$ with $x_k^{(j)} \rightarrow 0$, $l_k^{(j)} \rightarrow \infty$ and $1\leq m \leq \frac{T}{A_{\max}^{-n/2}\lambda_\gamma^\ast}$ such that 
\begin{equation}
v_k(x_k^{(j)}) = \max_{|x-x_k^{(j)}| \leq l_k^{(j)}\delta_k^{(j)}} v_k(x) \rightarrow +\infty 
 \label{eqn:lem:3y-1}
\end{equation}
for any $0\leq j \leq m-1$,  
\begin{equation}
B_{2l_k^{(i)}\delta_k^{(i)}}(x_k^{(i)}) \cap B_{2l_k^{(j)}\delta_k^{(j)}}(x_k^{(j)}) = \emptyset 
 \label{eqn:lem:3y-2}
\end{equation}
for any $k$ and $0 \leq i,j \leq m-1$ satisfying $i \neq j$,  
\begin{equation}
\left. \frac{\partial}{\partial t} v_k(ty+x_k^{(j)}) \right|_{t=1} < 0 
 \label{eqn:lem:3y-3}
\end{equation}
for any $k$, $0\leq j \leq m-1$ and $y$ satisfying $2r_\gamma^\ast\delta_k^{(j)} \leq |y| \leq 2l_k^{(j)}\delta_k^{(j)}$, 
\begin{align}
&\lim_{k\rightarrow\infty} \int_{B_{2l_k^{(j)}\delta_k^{(j)}}(x_k^{(j)})} (v_k)_+^{\frac{n(\gamma-1)}{2}} dx \nonumber\\
&= \lim_{k\rightarrow\infty} \int_{B_{l_k^{(j)}\delta_k^{(j)}}(x_k^{(j)})} (v_k)_+^{\frac{n(\gamma-1)}{2}} dx = A_0^{-n/2}\lambda_\gamma^\ast 
 \label{eqn:lem:3y-4}
\end{align}
for any $0\leq j \leq m-1$, and 
\begin{equation}
\sup_k \max_{x\in\overline{B_R}} \left\{ v_k(x) \min_{0 \leq j \leq m-1} |x-x_k^{(j)}|^q \right\} <+\infty, 
 \label{eqn:lem:3y-5}
\end{equation}
where $(\delta_k^{(j)})^{-q} = v_k(x_k^{(j)})$, $q=\frac{2}{\gamma -1}$, $A_0=A(0)$, and $r_\gamma^\ast$ is as in Theorem \ref{thm:class}. 
\end{lem}

\begin{lem}\label{lem:3z}
Assume the assumptions of Lemma \ref{lem:3y} and that there exist $\{ x_k^{(j)} \}_{j=0}^{m-1}$ and $\{ r_k^{(j)} \}_{j=0}^{m-1}$, $m\geq 1$, $r_k^{(j)}>0$, such that 
\begin{equation}
v_k(x_k^{(j)}) \rightarrow +\infty 
 \label{eqn:lem:3z-1}
\end{equation}
for any $0\leq j \leq m-1$, 
\begin{equation}
\lim_{k\rightarrow \infty} \frac{r_k^{(j)}}{\delta_k^{(j)}} = +\infty
 \label{eqn:lem:3z-2}
\end{equation}
for any $0\leq j \leq m-1$,  
\begin{equation}
B_{2r_k^{(i)}}(x_k^{(i)}) \cap B_{2r_k^{(j)}}(x_k^{(j)}) = \emptyset 
 \label{eqn:lem:3z-3}
\end{equation}
for any $k$ and $0 \leq i,j \leq m-1$ satisfying $i \neq j$,  
\begin{equation}
\sup_k \max_{x\in\overline{B_R}\setminus \cup_{j=0}^{m-1} B_{r_k^{(j)}}(x_k^{(j)})} 
\left\{ v_k(x) \min_{0\leq j \leq m-1} |x-x_k^{(j)}|^q \right\} <+\infty, 
 \label{eqn:lem:3z-4}
\end{equation}
and 
\begin{equation}
\lim_{k\rightarrow\infty} \int_{B_{2r_k^{(j)}}(x_k^{(j)})} (v_k)_+^{\frac{n(\gamma-1)}{2}} dx 
= \lim_{k\rightarrow\infty}\int_{B_{r_k^{(j)}}(x_k^{(j)})} (v_k)_+^{\frac{n(\gamma-1)}{2}} dx 
= \beta_j 
 \label{eqn:lem:3z-5}
\end{equation}
for some $\beta_j >0$, $0\leq j \leq m-1$. 
Then it holds that 
\begin{equation}
\lim_{k\rightarrow\infty} \int_{B_R} (v_k)_+^{\frac{n(\gamma-1)}{2}} dx = \sum_{j=0}^{m-1} \beta_j.  
 \label{eqn:lem:3z-6}
\end{equation}
\end{lem} 

The remainder of this section is devoted to the proof of Lemmas \ref{lem:3y}-\ref{lem:3z}. \\

\underline{{\it Proof of Lemma \ref{lem:3y}}}\ \ The proof consists of seven steps. \\ 

{\it Step 1.} We define $x_k^{(0)}$, $\delta_k^{(0)}$, $\tilde{v}_k^{(0)}$ and $\tilde{A}_k^{(0)}$ by 
\begin{align*}
&v_k(x_k^{(0)})=(\delta_k^{(0)})^{-q} = \max_{\overline{B_R}} v_k, \\ 
&\tilde{v}_k^{(0)}(x)=(\delta_k^{(0)})^q v_k(\delta_k^{(0)}x+x_k^{(0)}), \quad 
\tilde{A}_k^{(0)}(x)=A_k(\delta_k^{(0)}x+x_k^{(0)}), 
\end{align*}
where $q=\frac{2}{\gamma-1}$. 
From the scaling invariance, (\ref{eqn:lem:3x-1}) and (\ref{eqn:lem:3y-0}), it follows that  
\[
\begin{cases}
-\Delta\tilde{v}_k=\tilde{A}_k^{(0)}(x)(\tilde{v}_k)_+^\gamma,\quad \tilde{v}_{k} \leq 1\quad \mbox{in $B_{\frac{R}{4\delta_k^{(0)}}}$} \\ 
\int_{B_{\frac{R}{4\delta_k^{(0)}}}} (\tilde{v}_k)_+^{\frac{n(\gamma-1)}{2}} dx=\int_{B_{R/4}} (v_k)_+^{\frac{n(\gamma-1)}{2}} dx \leq T \\ 
\tilde{v}_k(0)=1. 
\end{cases}
\]
Using the elliptic regularity, Lemma \ref{lem:pre-2} and $\delta_k^{(0)} \rightarrow 0$, 
we obtain $\tilde{v} \in C^2({\bf R}^n)$, passing to a subsequence, such that $\tilde{v}_k^{(0)} \rightarrow \tilde{v}$ in $C_{loc}^{1+\alpha}({\bf R}^n)$ ($\alpha\in (0,1)$) and 
\[
\begin{cases}
-\Delta\tilde{v}=A_0\tilde{v}_+^\gamma,\quad \tilde{v}\leq 1\quad \mbox{in ${\bf R}^n$} \\ 
\int_{{\bf R}^n} \tilde{v}_+^{\frac{n(\gamma-1)}{2}} dx \leq T \\ 
\tilde{v}(0)=1, 
\end{cases}
\]
where $A_0=A(0)$. 
Since $\tilde{z}(x)=A_0^{\frac{1}{\gamma-1}}\tilde{v}(x)$ satisfies
\[
\begin{cases}
-\Delta\tilde{z}=\tilde{z}_+^\gamma,\quad \tilde{z}\leq A_0^{\frac{1}{\gamma-1}}\quad \mbox{in ${\bf R}^n$} \\ 
\int_{{\bf R}^n} \tilde{z}_+^{\frac{n(\gamma-1)}{2}} dx\leq A_0^{n/2}T\\ 
\tilde{z}(0)=A_0^{\frac{1}{\gamma-1}},
\end{cases}
\]
we obtain $\{ l_k^{(0)} \} \subset {\bf N}$, by virtue of Theorem \ref{thm:class} and a diagonal argument, such that $l_k^{(0)} \rightarrow \infty$ and
\begin{align*}
&\| \tilde{v}_k^{(0)} - \tilde{v} \|_{C^2(B_{2l_k^{(0)}})} \rightarrow 0 \\
&\int_{B_{2l_k^{(0)}\delta_k^{(0)}}(x_k^{(0)})} (v_k)_+^{\frac{n(\gamma-1)}{2}} dx 
= \int_{B_{2l_k^{(0)}}} (\tilde{v}_k^{(0)})_+^{\frac{n(\gamma-1)}{2}} dx \rightarrow A_0^{-n/2}\lambda_\gamma^\ast \\
&\int_{B_{l_k^{(0)}\delta_k^{(0)}}(x_k^{(0)})} (v_k)_+^{\frac{n(\gamma-1)}{2}} dx 
= \int_{B_{l_k^{(0)}}} (\tilde{v}_k^{(0)})_+^{\frac{n(\gamma-1)}{2}} dx \rightarrow A_0^{-n/2}\lambda_\gamma^\ast \\
&\left. \frac{\partial}{\partial t} v_k(ty+x_k^{(0)}) \right|_{t=1} < 0 
\quad \mbox{for any $y$ satisfying $2r_\gamma^\ast \delta_k^{(0)} \leq |y| \leq 2l_k^{(0)}\delta_k^{(0)}$.}
\end{align*} 

{\it Step 2.} It is clear that $\{x_k^{(0)} \}$ and $\{ l_k^{(0)} \}$ satisfy (\ref{eqn:lem:3y-1}) and (\ref{eqn:lem:3y-3})-(\ref{eqn:lem:3y-4}) for $m=1$. 

In turn, assume that $N$-sequences $\{ x_k^{(j)} \}_{j=0}^{N-1}$ and $\{ l_k^{(j)} \}_{j=0}^{N-1}$ satisfy (\ref{eqn:lem:3y-1})-(\ref{eqn:lem:3y-4}) for $m=N$,  
provided that (\ref{eqn:lem:3y-2}) is empty if $N=1$. 
Firstly, if 
\begin{equation}
\sup_k \max_{x\in \overline{B_R}} \{ v_k(x) \min_{0\leq j \leq N-1} |x-x_k^{(j)}|^q \} <+\infty
 \label{eqn:3y-6}
\end{equation}
then (\ref{eqn:lem:3y-1})-(\ref{eqn:lem:3y-5}) for $m=N$ hold, and so the lemma is true. 

Next we assume that (\ref{eqn:3y-6}) is false for the $N$-sequences above, which is supposed until Step 6 below is finished. 
Then there exists $\{ \overline{x}_k^{(N)} \}\subset B_R$ such that 
\begin{align}
P_k &\equiv v_k(\overline{x}_k^{(N)})\min_{0\leq j \leq N-1} |\overline{x}_k^{(N)} - x_k^{(j)}|^q \nonumber\\
&= \max_{x\in \overline{B_R}} \{ v_k(x) \min_{0\leq j \leq N-1} |x-x_k^{(j)}|^q \} \rightarrow +\infty.  
 \label{eqn:3y-7}
\end{align}
We put 
\begin{equation}
(\overline{\delta}_k^{(N)})^{-q} = v_k(\overline{x}_k^{(N)}) 
 \label{eqn:3y-9}
\end{equation}
and have
\begin{equation}
\frac{\min_{0\leq j \leq N-1} |\overline{x}_k^{(N)} - x_k^{(j)}|}{\overline{\delta}_k^{(N)}} = P_k^{1/q} \rightarrow +\infty 
 \label{eqn:3y-10}
\end{equation}
by (\ref{eqn:3y-7})-(\ref{eqn:3y-9}). 
Then it holds that 
\begin{align}
&\min_{0\leq j \leq N-1} |\overline{x}_k^{(N)} + \overline{\delta}_k^{(N)}x - x_k^{(j)}| \nonumber\\ 
&\geq -\overline{\delta}_k^{(N)}|x| + \min_{0\leq j \leq N-1} |\overline{x}_k^{(N)} - x_k^{(j)}| 
\geq \frac{1}{2}\min_{0\leq j \leq N-1} |\overline{x}_k^{(N)} - x_k^{(j)}| >0
 \label{eqn:3y-11}
\end{align}
for any $x$ satisfying 
\[
|x| \leq \frac{\min_{0\leq j \leq N-1} |\overline{x}_k^{(N)} - x_k^{(j)}|}{2\overline{\delta}_k^{(N)}} = \frac{P_k^{1/q}}{2} \equiv L_k\rightarrow +\infty. 
\]
Here we introduce 
\begin{equation}
\overline{v}_k(x)=(\overline{\delta}_k^{(N)})^q v_k(\overline{\delta}_k^{(N)}x+\overline{x}_k^{(N)}). 
 \label{eqn:3y-13}
\end{equation}
Note that 
\begin{align*}
v_k(\overline{\delta}_k^{(N)}x+\overline{x}_k^{(N)})
&\leq \frac{P_k}{\min_{0\leq j \leq N-1} |\overline{\delta}_k^{(N)}x+\overline{x}_k^{(N)}-x_k^{(j)}|^q} \\
&\leq \frac{2^q P_k}{\min_{0\leq j \leq N-1} |\overline{x}_k^{(N)}-x_k^{(j)}|^q}
=2^q(\overline{\delta}_k^{(N)})^{-q}
\end{align*}
for $x\in B_{L_k}$ by (\ref{eqn:3y-7}), (\ref{eqn:3y-11}) and (\ref{eqn:3y-10}). 
Thus, $\overline{v}_k=\overline{v}_k(x)$ satisfies 
\begin{equation}
\begin{cases}
-\Delta\overline{v}_k=\overline{A}_k(x)(\overline{v}_k)_+^\gamma,\quad \overline{v}_{k} \leq 2^q\quad \mbox{in $B_{L_k}$} \\ 
\int_{B_{L_k}} (\overline{v}_k)_+^{\frac{n(\gamma-1)}{2}} dx \leq T \\ 
\overline{v}_k(0) = 1,
\end{cases}
 \label{eqn:3y-14}
\end{equation}
where $\overline{A}_k(x)=A_k(\overline{\delta}_k^{(N)}x+\overline{x}_k^{(N)})$. 
Since $L_k\rightarrow+\infty$, passing to a subsequence, 
we obtain $\overline{v} \in C^2({\bf R}^n)$ such that $\overline{v}_k \rightarrow \overline{v}$ in $C_{loc}^{1+\alpha}({\bf R}^n)$ ($\alpha\in (0,1)$) and 
\[
\begin{cases}
-\Delta\overline{v}=A_0\overline{v}_+^\gamma,\quad \overline{v}\leq 2^q\quad \mbox{in ${\bf R}^n$} \\ 
\int_{{\bf R}^n} \overline{v}_+^{\frac{n(\gamma-1)}{2}} dx \leq T \\ 
\overline{v}(0)=1. 
\end{cases}
\]
Since $\overline{z}(x)=A_0^{\frac{1}{\gamma-1}}\overline{v}(x)$ satisfies 
\[
\begin{cases}
-\Delta\overline{z}=\overline{z}_+^\gamma,\quad \overline{z}\leq 2^q A_0^{\frac{1}{\gamma-1}}\quad \mbox{in ${\bf R}^n$} \\ 
\int_{{\bf R}^n} \overline{z}_+^{\frac{n(\gamma-1)}{2}} dx\leq A_0^{\frac{1}{\gamma-1}}T \\ 
\overline{z}(0)=A_0^{\frac{1}{\gamma-1}}, 
\end{cases}
\]
Theorem \ref{thm:class} yields 
\begin{align}
&\overline{v}(x)=\begin{cases}
A_0^{-\frac{1}{\gamma-1}}\mu^q\phi(\mu|x-\overline{x}_0|) & (|x-\overline{x}_0|\leq r_\gamma^\ast/\mu) \\
\frac{A_0^{-\frac{1}{\gamma-1}}\mu^{q-(n-2)}\lambda_\gamma^\ast}{\omega_{n-1}(n-2)}\left(\frac{1}{|x-\overline{x}_0|^{n-2}}-\frac{1}{(r_\gamma^\ast/\mu)^{n-2}}\right) & (|x-\overline{x}_0|\geq r_\gamma^\ast/\mu) 
\end{cases} 
 \label{eqn:3y-16} \\
&\int_{{\bf R}^n} \overline{v}_+^{\frac{n(\gamma-1)}{2}}dx=A_0^{-n/2}\lambda_\gamma^\ast
 \label{eqn:3y-16'}
\end{align}
for some 
\begin{equation}
\mu\in [\sqrt{A_0},2\sqrt{A_0}] \quad \mbox{and} \quad \overline{x}_0 \in B_{r_\gamma^\ast/\mu}.
 \label{eqn:3y-17}
\end{equation}

{\it Step 3.} Similarly to Step 1, there exists $\{ l_k^{(N)} \} \subset {\bf N}$ such that $l_k^{(N)} \rightarrow \infty$ and 
\begin{align}
&\| \overline{v}_k - \overline{v} \|_{C^2(B_{3l_k^{(N)}})} \rightarrow 0 
 \label{eqn:3y-18} \\ 
&\int_{B_{3l_k^{(N)}\overline{\delta}_k^{(N)}}(\overline{x}_k^{(N)})} (v_k)_+^{\frac{n(\gamma-1)}{2}} dx 
= \int_{B_{3l_k^{(N)}}} (\overline{v}_k)_+^{\frac{n(\gamma-1)}{2}} dx \rightarrow A_0^{-n/2}\lambda_\gamma^\ast
 \label{eqn:3y-19} \\ 
&\int_{B_{\frac{1}{4}l_k^{(N)}\overline{\delta}_k^{(N)}}(\overline{x}_k^{(N)})} (v_k)_+^{\frac{n(\gamma-1)}{2}} dx 
= \int_{B_{\frac{1}{4}l_k^{(N)}}} (\overline{v}_k)_+^{\frac{n(\gamma-1)}{2}} dx \rightarrow A_0^{-n/2}\lambda_\gamma^\ast
 \label{eqn:3y-20} \\ 
&\left. \frac{\partial}{\partial t} \overline{v}_k(ty+\overline{x}_0) \right|_{t=1} < 0 
\quad \mbox{for any $y$ satisfying $2r_\gamma^\ast \leq |y| \leq 3l_k^{(N)}$}. 
 \label{eqn:3y-21} 
\end{align}
Regarding (\ref{eqn:3y-21}) and $l_k^{(N)} \rightarrow \infty$, we take $y_k^{(N)} \in B_{2l_k^{(N)}}$ such that 
\begin{equation}
\overline{v}_k(\overline{x}_0+y_k^{(N)}) = \max_{y \in \overline{B_{3l_k^{(N)}}}} \overline{v}_k(\overline{x}_0+y), 
 \label{eqn:3y-22}
\end{equation}
and put 
\begin{equation}
x_k^{(N)} = \overline{x}_k^{(N)} + \overline{\delta}_k^{(N)}(\overline{x}_0+y_k^{(N)}). 
 \label{eqn:3y-23}
\end{equation}
Then it holds that 
\begin{align}
&y_k^{(N)} \rightarrow 0
 \label{eqn:3y-24} \\ 
&v_k(\overline{x}_k^{(N)}) \leq v_k(x_k^{(N)}) \leq 3^q v_k(\overline{x}_k^{(N)}). 
 \label{eqn:3y-25}
\end{align}
In fact, (\ref{eqn:3y-24}) follows from (\ref{eqn:3y-18}) and the fact that $\overline{v}=\overline{v}(x)$ attains its maximum at $x=\overline{x}_0$, recall (\ref{eqn:3y-16}). 
Also, (\ref{eqn:3y-25}) is derived from
\[
v_k(x_k^{(N)}) = (\overline{\delta}_k^{(N)})^{-q} \overline{v}_k(\overline{x}_0+y_k^{(N)}) 
\geq (\overline{\delta}_k^{(N)})^{-q} \overline{v}_k(0) 
= v_k(\overline{x}_k^{(N)}) 
\]
and 
\begin{align*}
v_k(x_k^{(N)}) = (\overline{\delta}_k^{(N)})^{-q} \overline{v}_k(\overline{x}_0+y_k^{(N)})  
= v_k(\overline{x}_k^{(N)}) \overline{v}_k(\overline{x}_0+y_k^{(N)}) 
\leq 3^q v_k(\overline{x}_k^{(N)}), 
\end{align*}
where we have used (\ref{eqn:3y-9}), (\ref{eqn:3y-13}), (\ref{eqn:3y-17})-(\ref{eqn:3y-18}), (\ref{eqn:3y-22})-(\ref{eqn:3y-23}) and $\overline{v}\leq 2^q$ in ${\bf R}^n$. \\ 

{\it Step 4. } We now claim 
\begin{align}
&\delta_k^{(N)} \leq \overline{\delta}_k^{(N)} \leq 3 \delta_k^{(N)}
 \label{eqn:3y-27} \\ 
&v_k(x_k^{(N)}) = \max_{|x-x_k^{(N)}| \leq l_k^{(N)}\delta_k^{(N)}} v_k(x) \rightarrow +\infty,  
 \label{eqn:3y-28}
\end{align}
where 
\[
\delta_k^{(N)} = (v_k(x_k^{(N)}))^{-1/q}.  
\]
Inequality (\ref{eqn:3y-27}) follows from (\ref{eqn:3y-25}). 
To show (\ref{eqn:3y-28}), we have only to prove the equality since the limit holds by (\ref{eqn:3y-25}) and $v_k(\overline{x}^{(N)})\rightarrow +\infty$. 
It holds that 
\begin{equation}
B_{l_k^{(N)}\delta_k^{(N)}}(x_k^{(N)}) \subset B_{3l_k^{(N)}\overline{\delta}_k^{(N)}}(\overline{x}_k^{(N)}+\overline{\delta}_k^{(N)}\overline{x}_0) 
 \label{eqn:3y-29}
\end{equation}
by (\ref{eqn:3y-27}) and  
\begin{align*}
|x_k^{(N)} - (\overline{x}_k^{(N)}+\overline{\delta}_k^{(N)}\overline{x}_0)| 
&\leq |x_k^{(N)}-\overline{x}_k^{(N)}| + \overline{\delta}_k^{(N)}|\overline{x}_0| 
= \overline{\delta}_k^{(N)} (|\overline{x}_0+y_k^{(N)}|+|\overline{x}_0|) \\ 
&\leq \overline{\delta}_k^{(N)} (3r_\gamma^\ast/\sqrt{A_0} + |y_k^{(N)}|) 
\leq \frac{3}{2} l_k^{(N)}\overline{\delta}_k^{(N)} 
\end{align*}
derived from (\ref{eqn:3y-23}), (\ref{eqn:3y-17}), (\ref{eqn:3y-24}) and $l_k^{(N)} \rightarrow \infty$. 
Using (\ref{eqn:3y-23}), (\ref{eqn:3y-22}), (\ref{eqn:3y-13}) and (\ref{eqn:3y-29}), we compute 
\begin{align*}
v_k(x_k^{(N)}) &= (\overline{\delta}_k^{(N)})^{-q} \overline{v}_k(\overline{x}_0+y_k^{(N)}) 
= \max_{y \in \overline{B_{3l_k^{(N)}}}} (\overline{\delta}_k^{(N)})^{-q} \overline{v}_k(\overline{x}_0+y) \\ 
&= \max_{y \in \overline{B_{3l_k^{(N)}\overline{\delta}_k^{(N)}}}} v_k(\overline{x}_k^{(N)}+\overline{\delta}_k^{(N)}\overline{x}_0+y) \\ 
&\geq \max_{x\in \overline{B_{l_k^{(N)}\delta_k^{(N)}}(x_k^{(N)})}} v_k(x) 
= \max_{|x-x_k^{(N)}| \leq l_k^{(N)}\delta_k^{(N)}} v_k(x) \geq v_k(x_k^{(N)}), 
\end{align*}
and hence the equality of (\ref{eqn:3y-28}) holds. \\ 

{\it Step 5.} Define 
\begin{equation}
\tilde{v}_k^{(N)}(x) = (\delta_k^{(N)})^q v_k(\delta_k^{(N)}x+x_k^{(N)}). 
 \label{eqn:3y-26}
\end{equation}
We next claim 
\begin{align}
&\| \tilde{v}_k^{(N)} - \tilde{v}^{(N)} \|_{C^2(B_{2l_k^{(N)}})} \rightarrow 0
 \label{eqn:3y-30} \\ 
&\int_{B_{2l_k^{(N)}\delta_k^{(N)}}(x_k^{(N)})} (v_k)_+^{\frac{n(\gamma-1)}{2}} dx 
= \int_{B_{2l_k^{(N)}}} (\tilde{v}_k^{(N)})_+^{\frac{n(\gamma-1)}{2}} dx \rightarrow A_0^{-n/2}\lambda_\gamma^\ast
 \label{eqn:3y-31} \\ 
&\int_{B_{l_k^{(N)}\delta_k^{(N)}}(x_k^{(N)})} (v_k)_+^{\frac{n(\gamma-1)}{2}} dx 
= \int_{B_{l_k^{(N)}}} (\tilde{v}_k^{(N)})_+^{\frac{n(\gamma-1)}{2}} dx \rightarrow A_0^{-n/2}\lambda_\gamma^\ast
 \label{eqn:3y-32} \\ 
&\left. \frac{\partial}{\partial t} v_k(ty+x_k^{(N)}) \right|_{t=1} < 0 
\quad \mbox{for any $y$ satisfying $2r_\gamma^\ast\delta_k^{(N)} \leq |y| \leq 2l_k^{(N)}\delta_k^{(N)}$}, 
 \label{eqn:3y-33}
\end{align}
where $\tilde{v}^{(N)}=\tilde{v}^{(N)}(x)$ is a function of the form (\ref{eqn:3y-16}). 
It is not difficult to check (\ref{eqn:3y-30}) and (\ref{eqn:3y-33}) similarly to Step 1. 
To prove (\ref{eqn:3y-31})-(\ref{eqn:3y-32}), it suffices to show that 
\begin{align}
&B_{2l_k^{(N)}\delta_k^{(N)}}(x_k^{(N)}) \subset B_{3l_k^{(N)}\overline{\delta}_k^{(N)}}(\overline{x}_k^{(N)})
 \label{eqn:3y-34} \\ 
&B_{\frac{1}{4}l_k^{(N)}\overline{\delta}_k^{(N)}}(\overline{x}_k^{(N)}) \subset B_{l_k^{(N)}\delta_k^{(N)}}(x_k^{(N)}) 
 \label{eqn:3y-35} 
\end{align}
by virtue of (\ref{eqn:3y-19})-(\ref{eqn:3y-20}). 
Relations (\ref{eqn:3y-27}), (\ref{eqn:3y-23})-(\ref{eqn:3y-24}) and (\ref{eqn:3y-17}) imply 
\begin{align*}
B_{2l_k^{(N)}\delta_k^{(N)}}(x_k^{(N)}) 
&\subset B_{2l_k^{(N)}\overline{\delta}_k^{(N)}}(x_k^{(N)}) \\
&\subset B_{2(l_k^{(N)}+2r_\gamma^\ast/\sqrt{A_0})\overline{\delta}_k^{(N)}}(\overline{x}_k^{(N)})
\subset B_{3l_k^{(N)}\overline{\delta}_k^{(N)}}(\overline{x}_k^{(N)}) 
\end{align*}
and 
\begin{align*}
B_{l_k^{(N)}\delta_k^{(N)}}(x_k^{(N)}) 
&\supset B_{\frac{1}{3}l_k^{(N)}\overline{\delta}_k^{(N)}}(x_k^{(N)}) \\
&\supset B_{(\frac{1}{3}l_k^{(N)}-2r_\gamma^\ast/\sqrt{A_0})\overline{\delta}_k^{(N)}}(\overline{x}_k^{(N)})
\supset B_{\frac{1}{4}l_k^{(N)}\overline{\delta}_k^{(N)}}(\overline{x}_k^{(N)}) 
\end{align*}
for $k \gg 1$, and then (\ref{eqn:3y-34})-(\ref{eqn:3y-35}) follow. \\ 

{\it Step 6.} We are now in a position to show that 
$\{ x_k^{(j)} \}_{j=0}^{N}$ and $\{ l_k^{(j)} \}_{j=0}^{N}$ satisfy the properties (\ref{eqn:lem:3y-1})-(\ref{eqn:lem:3y-4}) for $m=N+1$. 
Properties (\ref{eqn:lem:3y-1}) and (\ref{eqn:lem:3y-3})-(\ref{eqn:lem:3y-4}) 
follow from (\ref{eqn:3y-28}), (\ref{eqn:3y-30})-(\ref{eqn:3y-33}) and the hypothesis of induction. 
The proof of (\ref{eqn:lem:3y-2}) is reduced to showing that of 
\begin{equation}
|x_k^{(N)} - x_k^{(j)}| > 2r_\gamma^\ast (\delta_k^{(N)} + \delta_k^{(j)}) 
 \label{eqn:3y-36} 
\end{equation}
for any $0\leq j \leq N-1$. 
Indeed, if (\ref{eqn:3y-36}) and 
\[
B_{2l_{k_0}^{(N)}\delta_{k_0}^{(N)}}(x_{k_0}^{(N)}) \cap B_{2l_{k_0}^{(j)}\delta_{k_0}^{(j)}}(x_{k_0}^{(j)}) \neq \emptyset 
\]
occur simultaneously for some $k_0$ and $0\leq j \leq N-1$, then there is a point $z_0$ on the segment $\overline{x_{k_0}^{(N)}x_{k_0}^{(j)}}$, such that 
\[ 
|z_0-x_{k_0}^{(N)}|\geq 2r_\gamma^\ast\delta_{k_0}^{(N)} \quad \mbox{and} \quad |z_0-x_{k_0}^{(j)}|\geq 2r_\gamma^\ast\delta_{k_0}^{(j)}, 
\]
which is impossible because of (\ref{eqn:lem:3y-3}) for $m=N+1$. 
Therefore, (\ref{eqn:3y-36}) implies (\ref{eqn:lem:3y-2}) for $m=N+1$. 

To show (\ref{eqn:3y-36}) by contradiction, assume that there exists $\{x_{k'}^{(j)}\}\subset\{x_k^{(j)}\}$ such that 
\[
|x_{k'}^{(N)}-x_{k'}^{(j)}| \leq 2r_\gamma^\ast (\delta_{k'}^{(N)}+\delta_{k'}^{(j)}). 
\]
Noting that $\delta_{k'}^{(N)} = \delta_{k'}^{(j)}$ for $k'\gg 1$ by (\ref{eqn:3y-28}) and the hypothesis of induction, we calculate 
\begin{align*}
&|x_{k'}^{(N)}-x_{k'}^{(j)}| 
\geq |\overline{x}_{k'}^{(N)}-x_{k'}^{(j)}| - \overline{\delta}_{k'}^{(N)} |\overline{x}_0+y_{k'}^{(N)}| \\
&\geq (P_{k'}^{1/q}-|\overline{x}_0+y_{k'}^{(N)}|) \overline{\delta}_{k'}^{(N)} 
\geq \frac{1}{2}P_{k'}^{1/q}\overline{\delta}_{k'}^{(N)} 
\geq \frac{1}{2}P_{k'}^{1/q}\delta_{k'}^{(N)} 
\end{align*}
for $k' \gg 1$ by (\ref{eqn:3y-23}), (\ref{eqn:3y-10}), (\ref{eqn:3y-17}), (\ref{eqn:3y-24}) and (\ref{eqn:3y-27}). 
Consequently, $\frac{1}{2}P_{k'}^{1/q} \leq 4r_\gamma^\ast$ for $k' \gg 1$, which is false by (\ref{eqn:3y-10}). 
Hence (\ref{eqn:lem:3y-2}) is shown for $m=N+1$. \\ 

{\it Step 7.} We have shown 
that there exist $\{ x_k^{(j)} \}_{j=0}^{N}$ and $\{ l_k^{(j)} \}_{j=0}^{N}$ satisfying (\ref{eqn:lem:3y-1})-(\ref{eqn:lem:3y-4}) for $m=N+1$ 
under the assumptions that $\{ x_k^{(j)} \}_{j=0}^{N-1}$ and $\{ l_k^{(j)} \}_{j=0}^{N-1}$ satisfy the same properties for $m=N$ and that (\ref{eqn:3y-6}) is false.  
We can continue the procedure developed above as far as (\ref{eqn:3y-6}) fails. 
On the other hand, the procedure must end with finite times, $\left[\frac{T}{A_{\max}^{-n/2}\lambda_\gamma^\ast}\right]$ times at most. 
Eventually (\ref{eqn:lem:3y-5}) holds for some $m$, and the lemma is established. \qed \\ 

\underline{{\it Proof of Lemma \ref{lem:3z}}}\ \ We shall prove the lemma by induction. 
The proof is divided into four steps. \\

{\it Step 1.} It holds by (\ref{eqn:lem:3x-2}) and (\ref{eqn:lem:3z-1}) that 
\begin{equation}
x_k^{(j)} \rightarrow 0 \quad \mbox{for $0 \leq j \leq m-1$}. 
 \label{eqn:3z-2}
\end{equation}

First of all, we shall prove the lemma for $m=1$. 
Without loss of generality, we may assume 
\begin{equation}
x_k^{(0)} = 0 \quad \mbox{for any $k$},  
 \label{eqn:3z-3}
\end{equation}
besides $r_k^{(0)} \rightarrow 0$, for, if the latter is not the case then the lemma follows from (\ref{eqn:lem:3x-2}) and (\ref{eqn:lem:3z-5}). 
Since (\ref{eqn:lem:3z-4}) holds for $m=1$, estimate (\ref{eqn:lem:3-2-2}) gives 
\begin{equation*}
\sup_{\partial B_r} v_k \leq r^{-q} C_{\const\insa} - v_k(0) C_{\const\insb}  
= r^{-q} C_{\suba} - (\delta_k^{(0)})^{-q} C_{\subb}
 \label{eqn:3z-4} 
\end{equation*} 
for any $r\in [2r_k^{(0)},R/2]$, and for some $C_{\suba}>0$ and $C_{\subb}>0$ independent of $r$, where $q=\frac{2}{\gamma-1}$. 
Then we get   
\begin{equation}
v_k \leq 0 \quad \mbox{on $\partial B_r$} \quad \mbox{if $(\frac{r}{\delta_k^{(0)}})^q \geq \frac{C_{\suba}}{C_{\subb}}$} 
 \label{eqn:3z-5} 
\end{equation}
for $k \gg 1$. 
Relations (\ref{eqn:3z-5}), (\ref{eqn:lem:3z-2}), (\ref{eqn:lem:3z-5}) and (\ref{eqn:lem:3x-2}) admit (\ref{eqn:lem:3z-6}) for $m=1$. \\ 

{\it Step 2.} In the following, we shall prove the lemma for $m \geq 2$. 
We may assume (\ref{eqn:3z-3}) and 
\begin{equation}
d_k = |x_k^{(0)}-x_k^{(1)}| = \min_{0\leq i,j \leq m-1, \ i\neq j} |x_k^{(i)}-x_k^{(j)}| 
 \label{eqn:3z-6}
\end{equation}
by relabeling the indices. 
There are two possibilities: 

\begin{description}
\item{{\it Case 1.}} There exists $R_1\geq 1$ such that 
\begin{equation}
|x_k^{(i)}-x_k^{(j)}| \leq R_1 d_k 
 \label{eqn:3z-7}
\end{equation}
for any k and $i \neq j$, $0\leq i,j \leq m-1$. 
\item{{\it Case 2.}} There exist $J \subset \{ 0,1,\cdots,m-1 \}$ and $R_2 \geq 1$ such that 
\begin{align}
& \{0,1\} \subset J
 \label{eqn:3z-8} \\ 
& |x_k^{(j)}| \leq R_2 d_k \quad \mbox{for any $k$ and $j \in J$}
 \label{eqn:3z-9} \\ 
& \lim_{k \rightarrow \infty} \frac{|x_k^{(j)}|}{d_k} = +\infty \quad \mbox{for any $j \notin J$.} 
 \label{eqn:3z-10} 
\end{align}
\end{description} 

We shall deal with {\it Case 1} and {\it Case 2} in {\it Step 3} and {\it Step 4} below, respectively. \\ 

{\it Step 3.} To show the lemma for {\it Case 1}, it suffices to say that 
\begin{equation}
\lim_{k\rightarrow \infty} \int_{B_{4R_1 d_k}} (v_k)_+^{\frac{n(\gamma-1)}{2}} dx 
= \lim_{k\rightarrow \infty} \int_{B_{2R_1 d_k}} (v_k)_+^{\frac{n(\gamma-1)}{2}} dx 
= \sum_{j=0}^{m-1}\beta_j. 
 \label{eqn:3z-11}
\end{equation}
To see this sufficiency, we put 
\begin{equation}
\begin{cases}
{x'}_k^{(0)} = x_k^{(0)} =0, \quad {r'}_k^{(0)} = 2R_1 d_k, \\ 
\beta_0' = \sum_{j=0}^{m-1} \beta_j, \quad {\delta'}_k^{(0)} = \delta_k^{(0)}, \quad m'=1, 
\end{cases}
 \label{eqn:3z-12}
\end{equation} 
and check that the lemma holds for $m=m'=1$. 
By the hypothesis of induction, we have only to confirm (\ref{eqn:lem:3z-1})-(\ref{eqn:lem:3z-2}) and (\ref{eqn:lem:3z-4})-(\ref{eqn:lem:3z-5}) for the quantities in (\ref{eqn:3z-12}). 
It is clear that (\ref{eqn:lem:3z-5}) is equivalent to (\ref{eqn:3z-11}). 
Also (\ref{eqn:lem:3z-1})-(\ref{eqn:lem:3z-2}) and (\ref{eqn:lem:3z-4}) follow from 
\begin{align*}
&v_k({x'}_k^{(0)}) = v_k(x_k^{(0)}) \rightarrow +\infty \\
&\frac{{r'}_k^{(0)}}{{\delta'}_k^{(0)}} = \frac{2R_1 |x_k^{(1)}|}{\delta_k^{(0)}} \geq \frac{2R_1 r_k^{(0)}}{\delta_k^{(0)}} \rightarrow +\infty \\ 
&\cup_{j=0}^{m-1} B_{r_k^{(j)}}(x_k^{(j)}) \subset B_{R_1 d_k+\max_{0\leq j \leq m-1}r_k^{(j)}} \subset B_{2R_1 d_k} = B_{{r'}_k^{(0)}}({x'}_k^{(0)}), 
\end{align*}
derived from (\ref{eqn:3z-6})-(\ref{eqn:3z-7}), (\ref{eqn:3z-12}) and the hypothesis of induction. 
Hence (\ref{eqn:3z-11}) implies the lemma. 

Now we shall show (\ref{eqn:3z-11}). 
we introduce 
\begin{equation}
\begin{cases} 
&\tilde{v}_k(x) = d_k^q v_k(d_k x) \quad \mbox{for $|x| \leq \frac{R}{d_k}$} \\ 
&\tilde{x}_k^{(j)} = x_k^{(j)} / d_k \quad \mbox{for $0 \leq j \leq m-1$} \\ 
&(\tilde{\delta}_k^{(j)})^{-q} = \tilde{v}_k(\tilde{x}_k^{(j)}) = d_k^q v_k(x_k^{(j)}) = ( \frac{d_k}{\delta_k^{(j)}} )^q \quad \mbox{for $0 \leq j \leq m-1$} \\ 
&\tilde{r}_k^{(j)} = r_k^{(j)} /d_k \quad \mbox{for $0 \leq j \leq m-1$}, 
\end{cases}
 \label{eqn:3z-13} 
\end{equation}
and use (\ref{eqn:lem:3z-1})-(\ref{eqn:lem:3z-5}) of the original sequences to find 
\begin{equation}
\begin{cases}
&\tilde{x}_k^{(0)} = 0 \\ 
&\tilde{v}_k(\tilde{x}_k^{(j)}) \rightarrow +\infty \quad \mbox{for $0 \leq j \leq m-1$} \\ 
&\frac{\tilde{r}_k^{(j)}}{\tilde{\delta}_k^{(j)}} \rightarrow +\infty \quad \mbox{for $0 \leq j \leq m-1$} \\ 
&B_{\tilde{r}_k^{(i)}}(\tilde{x}_k^{(i)}) \cap B_{\tilde{r}_k^{(j)}}(\tilde{x}_k^{(j)}) = \emptyset \quad \mbox{for $i \neq j$, $0 \leq j \leq m-1$}, \\ 
&\sup_k \max_{x\in\overline{B_{R/d_k}} \setminus \cup_{j=0}^{m-1} B_{\tilde{r}_k^{(j)}}(\tilde{x}_k^{(j)})} 
\{ \tilde{v}_k(x) \min_{0 \leq j \leq m-1} |x-\tilde{x}_k^{(j)}|^q \} <+\infty \\
&\lim_{k\rightarrow \infty} \int_{B_{2\tilde{r}_k^{(j)}}(\tilde{x}_k^{(j)})} (\tilde{v}_k)_+^{\frac{n(\gamma -1)}{2}} dx 
= \lim_{k\rightarrow \infty} \int_{B_{\tilde{r}_k^{(j)}}(\tilde{x}_k^{(j)})} (\tilde{v}_k)_+^{\frac{n(\gamma -1)}{2}} dx \\ 
&\quad\quad\quad\quad\quad\quad\quad \ \ \! 
= \beta_j \quad \mbox{for $0 \leq j \leq m-1$}. 
\end{cases}
 \label{eqn:3z-14}
\end{equation}
We may assume 
\begin{equation} 
\tilde{x}_k^{(j)} \rightarrow \tilde{x}^{(j)} \in \overline{B_{R_1}} 
 \label{eqn:3z-15} 
\end{equation} 
for any $0 \leq j \leq m-1$ by virtue of (\ref{eqn:3z-7}) and (\ref{eqn:3z-13}). 
Then the scaling invariance of (\ref{eqn:3z-13}), (\ref{eqn:3z-15}) and the second property of (\ref{eqn:3z-14}) allow us to apply Proposition \ref{pro:bm}, 
so that 
\begin{equation}
\tilde{v}_k \rightarrow -\infty \quad \mbox{locally uniformly in ${\bf R}^n \setminus \{ \tilde{x}^{(0)}, \cdots, \tilde{x}^{(m-1)} \}$}.
 \label{eqn:3z-16}
\end{equation}
It follows from (\ref{eqn:3z-6})-(\ref{eqn:3z-7}) and (\ref{eqn:3z-13}) that 
\begin{equation}
1 \leq |\tilde{x}_k^{(i)}-\tilde{x}_k^{(j)}| \leq R_1 \quad \mbox{for $i \neq j$, $0 \leq j \leq m-1$}. 
 \label{eqn:3z-17}
\end{equation} 
For each $j$, we have either
\[
\mbox{(i) $\limsup_{k\rightarrow \infty} \tilde{r}_k^{(j)} > 0$ \ \ or\ \ (ii) $\lim_{k\rightarrow \infty} \tilde{r}_k^{(j)} = 0$}. 
\]
If (i) occurs then  
\begin{equation}
\int_{B(\tilde{x}^{(j)},1/2)} (\tilde{v}_k)_+^{\frac{n(\gamma -1)}{2}} dx \rightarrow \beta_j 
 \label{eqn:3z-18}
\end{equation}
for each $j$ by (\ref{eqn:3z-16})-(\ref{eqn:3z-17}). 
If (ii) occurs then for  
\[
{x'}_k^{(0)} = \tilde{x}_k^{(j)} =0, \quad {r'}_k^{(0)} = \tilde{r}_k^{(j)}, \quad \beta_0' = \beta_j, \quad {\delta'}_k^{(0)} = \tilde{\delta}_k^{(j)}, \quad m'=1, 
\]
the assumptions of the lemma hold for $m=m'=1$, which implies (\ref{eqn:3z-18}) by {\it Step 1} and (\ref{eqn:3z-16})-(\ref{eqn:3z-17}). 
Consequently, (\ref{eqn:3z-16})-(\ref{eqn:3z-18}) hold for both cases, and thus (\ref{eqn:3z-11}) is shown. \\ 

{\it Step 4.} Without loss of generality, we may assume 
\[
J=\{0,1,\cdots,l-1\} \quad \mbox{for some $2\leq l \leq m-1$}. 
\]
Similarly to {\it Step 3}, we obtain 
\begin{equation}
\lim_{k\rightarrow \infty} \int_{B_{4R_2 d_k}} (v_k)_+^{\frac{n(\gamma-1)}{2}} dx 
= \lim_{k\rightarrow \infty} \int_{B_{2R_2 d_k}} (v_k)_+^{\frac{n(\gamma-1)}{2}} dx 
= \sum_{j=0}^{l-1}\beta_j 
 \label{eqn:3z-19}
\end{equation}
since it holds that   
\begin{align*}
&B_{4R_2 d_k}(x_k^{(j)}) \subset B_{5R_2 d_k} \quad \mbox{for any $k$ and $j \in J$} \\ 
&|x_k^{(j)}| \geq 6R_2 d_k \quad \mbox{for $k \gg 1$ and $j \notin J$} 
\end{align*}
by (\ref{eqn:3z-9}) and (\ref{eqn:3z-10}). 
We put 
\begin{align*}
&\tilde{v}_k(x) = d_k^q v_k(d_k x) \quad \mbox{for $|x| \leq 5R_2$} \\ 
&{x'}_k^{(0)} = x_k^{(0)} =0, \quad {r'}_k^{(0)} = 2R_2d_k, \quad \beta_0' = \sum_{j=0}^{l-1} \beta_j, \quad {\delta'}_k^{(0)} = \delta_k^{(0)},  
\end{align*}
and find by the hypothesis of induction that for 
\begin{equation}
\begin{cases}
&\{ {x'}_k^{(0)} \} \cup \{ x_k^{(j)} \}_{j=l}^{m-1}, \quad \{ {r'}_k^{(0)} \} \cup \{ r_k^{(j)} \}_{j=l}^{m-1}, \\ 
&\beta_0',\beta_l,\cdots,\beta_{m-1}, \quad {\delta'}_k^{(0)},\delta_k^{(l)},\cdots,\delta_k^{(m-1)},  
\end{cases}
 \label{eqn:3z-20} 
\end{equation}
if the assumptions of the lemma are satisfied for $m-l+1$ then 
\[
\int_{B_R} (v_k)_+^{\frac{n(\gamma-1)}{2}} dx \rightarrow \beta_0' + \sum_{j=l}^{m-1} \beta_j = \sum_{j=0}^{m-1} \beta_j. 
\]

Therefore, in order to show the lemma for {\it Case 2}, 
we have only to show that the quantities in (\ref{eqn:3z-20}) satisfy (\ref{eqn:lem:3z-1})-(\ref{eqn:lem:3z-5}). 
We shall show (\ref{eqn:lem:3z-3}) here since it is not difficult to check (\ref{eqn:lem:3z-1})-(\ref{eqn:lem:3z-2}) and (\ref{eqn:lem:3z-4})-(\ref{eqn:lem:3z-5}). 
To this end, assume the contrary, that is, 
\begin{equation}
B_{2R_2 d_k} \cap B_{r_k^{(j_0)}}(x_k^{(j_0)}) \neq \emptyset 
 \label{eqn:3z-21} 
\end{equation}
for some $l \leq j_0 \leq m-1$ and for $k \gg 1$. 
It holds that 
\[
\limsup_{k\rightarrow \infty} \frac{r_k^{(j_0)}}{|x_k^{(j_0)}|} \leq 1 
\]
by ${x'}_k^{(0)} = x_k^{(0)} =0 \notin B(x_k^{(j_0)},r_k^{(j_0)})$. 
We have, on the other hand, 
\[
2R_2 d_k + r_k^{(j_0)} \geq |x_k^{(j_0)}| 
\]
for $k \gg 1$ by (\ref{eqn:3z-21}).  
This inequality and (\ref{eqn:3z-10}) imply  
\[
\liminf_{k\rightarrow \infty} \frac{r_k^{(j_0)}}{|x_k^{(j_0)}|} \geq 1,   
\]
and thus 
\begin{equation}
\lim_{k\rightarrow \infty} \frac{r_k^{(j_0)}}{|x_k^{(j_0)}|} = 1. 
 \label{eqn:3z-22}
\end{equation}
Furthermore, it holds that 
\begin{equation}
\frac{r_k^{(j_0)}}{d_k} = \frac{r_k^{(j_0)}}{|x_k^{(j_0)}|} \cdot \frac{|x_k^{(j_0)}|}{d_k} \rightarrow +\infty 
 \label{eqn:3z-23}
\end{equation}
by (\ref{eqn:3z-10}) and (\ref{eqn:3z-22}). 
We organize (\ref{eqn:3z-9}) and (\ref{eqn:3z-22})-(\ref{eqn:3z-23}) to obtain 
\begin{equation*}
B_{r_k^{(0)}} \subset B_{2R_2 d_k} \subset B_{2r_k^{(j_0)}}(x_k^{(j_0)})
\end{equation*}
for $k\gg 1$, which implies 
\begin{align*}
\beta_{j_0} &= \lim_{k\rightarrow \infty} \int_{B_{2r_k^{(j_0)}}(x_k^{(j_0)})} (v_k)_+^{\frac{n(\gamma-1)}{2}} dx \\ 
&\geq \lim_{k\rightarrow \infty} \int_{B_{r_k^{(0)}}} (v_k)_+^{\frac{n(\gamma-1)}{2}} dx 
+ \int_{B_{r_k^{(j_0)}}(x_k^{(j_0)})} (v_k)_+^{\frac{n(\gamma-1)}{2}} dx 
= \beta_0 + \beta_{j_0} > \beta_{j_0}, 
\end{align*}
a contradiction. 
Hence we obtain (\ref{eqn:lem:3z-3}) for the quantities in (\ref{eqn:3z-20}). 
The proof is complete. \qed


\end{document}